\input amstex
\documentstyle{amsppt} 
\magnification=\magstep1
\leftheadtext{YOSHIAKI FUKUMA}
\rightheadtext{A LOWER BOUND FOR $K_{X}L$ OF 
QUASI-POLARIZED SURFACES}

\topmatter 
\title  A LOWER BOUND FOR $\bold{K_{X}L}$ OF QUASI-POLARIZED SURFACES
$\bold{(X,L)}$ WITH NON-NEGATIVE KODAIRA DIMENSION
\endtitle
\author YOSHIAKI FUKUMA
\endauthor
\address Department of Mathematics,
 Faculty of Science, 
 Tokyo Institute of Technology, 
 Oh-okayama, Meguro-ku, Tokyo 152,
 Japan
\endaddress
\email fukuma\@math.titech.ac.jp
\endemail
\date {July, 1996}\enddate
\abstract
Let $X$ be a smooth projective surface over the complex number field
and let $L$ be a nef-big divisor on $X$.
Here we consider the following conjecture;
If the Kodaira dimension $\kappa(X)\geq 0$, then $K_{X}L\geq 2q(X)-4$,
where $q(X)$ is the irregularity of $X$.
In this paper, we prove that this conjecture is true if
(1) the case in which $\kappa(X)=0$ or 1, (2) the case in which
$\kappa(X)=2$ and $h^{0}(L)\geq 2$, or
(3) the case in which $\kappa(X)=2$, $X$ is minimal, $h^{0}(L)=1$, and
$L$ satisfies some conditions.
\endabstract
\endtopmatter

\document
\subhead{\bf $\S$ 0. Introduction}
\endsubhead\par

Let $X$ be a smooth projective manifold over $\Bbb{C}$ with 
$\operatorname{dim}X\geq 2$,
and $L$ a Cartier divisor on $X$.
Then $(X,L)$ is called a pre-polarized manifold.
In particular, if $L$ is ample (resp. nef-big),
then $(X,L)$ is said to be a polarized (resp. quasi-polarized)
manifold.
We define the sectional genus $g(L)$ of a pre-polarized manifold
$(X,L)$ is defined by the following formula;
$$g(L)=1+\frac{1}{2}(K_{X}+(n-1)L)L^{n-1},$$
where $K_{X}$ is the canonical divisor of $X$.
\flushpar
Then this is the following conjecture.
\proclaim{Conjecture 0}
Let $(X,L)$ be a quasi-polarized manifold.
Then $g(L)\geq q(X)$, where $q(X)=\operatorname{dim}H^{1}(X,\Cal{O}_{X})$.
\endproclaim
\flushpar
In this paper, we consider the case in which $X$ is a smooth surface.
If $\dim X=2$ and $h^{0}(L)>0$, then this Conjecture is true.
But in general, it is unknown whether this Conjecture is true or
not.
In the papers \cite{Fk1} and \cite{Fk4}, 
the author proved that $L^{2}\leq 4$ if $L$ is ample, $g(L)=q(X)$,
$h^{0}(L)>0$ and $\kappa(X)\geq 0$.
By this result, we think that the degree of $(X,L)$ is bounded
from above by using $m=g(L)-q(X)$ if $\kappa(X)\geq 0$.
By studying some examples of $(X,L)$, we conjectured the following.
\proclaim{Conjecture 1}
If $(X,L)$ is a quasi-polarized surface with $\kappa(X)\geq 0$.
\flushpar
Then $L^{2}\leq 2m+2$ if $g(L)=q(X)+m$.
\endproclaim
\flushpar
We remark that $m$ is non-negative integer if $h^{0}(L)>0$.
This Conjecture is equivalent to the following Conjecture.
\proclaim{Conjecture 1$'$} 
If $(X,L)$ is a quasi-polarized surface with $\kappa(X)\geq 0$.
\flushpar
Then $K_{X}L\geq 2q(X)-4$.
\endproclaim
This Conjecture 1$'$ is thought to be an generalization of the fact that 
$\operatorname{deg}K_{C}=2g(C)-2$ if $C$ is a smooth projective curve.
\flushpar
In this paper, we consider the above Conjecture.
Main results are the following.
\proclaim{Main Theorem 1}
Let $(X,L)$ be a quasi-polarized surface with $\kappa(X)=0$ or 1.
Then $K_{X}L\geq 2q(X)-4$.
\flushpar
If this equality holds and $(X,L)$ is $L$-minimal, then
$(X,L)$ is one of the following;
\roster
\item $\kappa(X)=0$ case.\flushpar
$X$ is an Abelian surface and $L$ is any nef and big divisor.
\item $\kappa(X)=1$ case. \flushpar
$X\cong F\times C$ and $L\equiv C+(m+1)F$,
where $F$ and $C$ are smooth curves with $g(C)\geq 2$ and $g(F)=1$,
and $m=g(L)-q(X)$.
\endroster
\endproclaim
(See Theorem 2.1.)
\proclaim{Main Theorem 2}
Let $(X,L)$ be a quasi-polarized surface with $\kappa(X)=2$ and
$h^{0}(L)\geq 2$.
Then $K_{X}L\geq 2q(X)-2$.
\flushpar
If this equality holds and $(X,L)$ is $L$-minimal,
then $(X,L)$ is the following;
\flushpar
$X\cong F\times C$ and $L\equiv C+2F$,
where $F$ and $C$ are smooth curves with $g(F)=2$ and $g(C)\geq 2$.
\endproclaim
(See Theorem 3.1)
\proclaim{Main Theorem 3}
Let $X$ be a smooth surface of general type
and let $D$ be a nef-big effective divisor with $h^{0}(D)=1$
on $X$.
If $D$ is not the following type ($\star$), then
$K_{X}D\geq 2q(X)-4$;
\roster
\item"($\star$)"
$D=C_{1}+\sum\limits_{j\geq 2}r_{j}C_{j}$ ;
$C_{1}^{2}>0$ and the intersection matrix 
$\Vert C_{j}, C_{k}\Vert_{j\geq 2, k\geq 2}$
of $\sum\limits_{j\geq 2}r_{j}C_{j}$ is
negative semidefinite.
\endroster
\endproclaim
(See $\S$ 4.)
\proclaim{Main Theorem 4}
Let $X$ be a minimal smooth projective surface with
$\kappa(X)=2$ and let $D$ be a nef-big effective divisor on $X$
such that $D$ is the type ($\star$).
Then $D^{2}\leq 4m+4$ if $m=g(D)-q(X)$.
\endproclaim
We remark that the classification of polarized surfaces $(X,L)$ with
$\kappa(X)\geq 1$ and $K_{X}L\leq 2$ is obtained by \cite{DP}.
\flushpar
{\bf Acknowledgement}\flushpar
The author would like to express his hearty gratitude
to Professor Takao Fujita for giving some useful advice and comments.

\subhead{$\S$ 1. Preliminaries}
\endsubhead
\definition{Definition 1.1}(\cite{Fk1}.) 
Let $(X,L)$ be a quasi-polarized surface.
\roster
\item We call $(X_{1},L_{1})$ a minimalization of $(X,L)$ 
if $\varphi : X\to X_{1}$ is a minimal model of $X$ 
and $L_{1}=\varphi_{*}L$ in the sense of cycle theory. 
(We remark that $L_{1}$ is nef and big (resp. ample) on $X_{1}$ if so is $L$.)
\item We say that $(X,L)$ is $L$-minimal if $L.E>0$ for any (-1)-curve $E$
on $X$. For any quasi-polarized surface $(X,L)$, there exists 
a birational morphism $\rho:(X,L) \to (X_{0},L_{0})$ such that 
$L=\rho^{*}L_{0}$ and $(X_{0},L_{0})$
is $L_{0}$-minimal. Then we call $(X_{0},L_{0})$ an $L$-minimalization
of $(X,L)$.
\endroster
\enddefinition
\proclaim{Lemma 1.2 (Debarre)}
Let $X$ be a minimal surface of general type with $q(X)\geq 1$.
Then $K_{X}^{2}\geq 2p_{g}(X)$.
(Hence $K_{X}^{2}\geq 2q(X)$ for any minimal surface of general type.)
\endproclaim
\demo{Proof}
See \cite{D}. \qed
\enddemo
\proclaim{Theorem 1.3 (\cite{Fk3})}
Let $(X,L)$ be an $L$-minimal quasi-polarized surface with $\kappa(X)\geq 0$.
If $h^{0}(L)\geq 2$, then $(X,L)$ satisfies one of the following conditions.
\roster
\item  $g(L)\geq 2q(X)-1$.
\item  For any linear pencil $\Lambda\subseteq |L|$,
the fixed part $Z(\Lambda)$ of $\Lambda$ is not zero 
and $\operatorname{Bs}\Lambda_{M}=\phi$, where $\Lambda_{M}$ is
movable part of $\Lambda$.
Let $f:X\to C$ be the fiber space induced by $\Lambda_{M}$.
Then $g(L)\geq g(C)+2g(F)\geq q(X)+g(F)$,
$g(C)\geq 2$, $LF=1$ and $L-aF$ is numerically equivalent
to an effective divisor for a general fiber $F$ of $f$,
where $a\geq 2$.
\endroster
\endproclaim
\demo{Proof}
See Theorem 3.1 in \cite{Fk3}. \qed
\enddemo
\proclaim{Lemma 1.4}
Let $f:X\to C$ be a relatively minimal elliptic fibration with $q(X)=g(C)+1$.
If $LF=1$ for a nef-big divisor $L$ on $X$, 
then $X\cong F\times C$ and $f:X\to C$ is the natural projection, 
where $F$ is a general fiber of $f$.
\endproclaim
\demo{Proof} (See \cite{Fj3})
By hypothesis $f$ is a quasi-bundle 
(see Lemma 1.5 and Lemma 1.6 in \cite{S}).
Let $\Sigma\subset C$ be the singular locus of $f$ and $U=C-\Sigma$.
We fix an elliptic curve $E\cong f^{-1}(x)$ for $x\in U$.
Then by \cite{Fj3}, we have a map 
$\varphi:\pi_{1}(U)\to \operatorname{Aut}(E,C'_{E})$.
Since the translations of $E$ preserving $L_{E}$ are 
of order $d=\operatorname{deg}L_{E}$ by Abel's Theorem,
$\operatorname{Aut}(E,L_{E})$ is finite group.
Let $G=\operatorname{Im}\varphi$.
Then by \cite{Fj3}, there exists a Galois covering
$\pi: D\to C$ such that $G=\operatorname{Gal}(D/C)$
acts effectively on the polarized pair $(E,L_{E})$
and $X\cong (D\times E)/G$,
where $D$ is a smooth projective curve.
Since $q(X)=g(C)+1$, we have $g(E/G)=1$. 
Hence $G$ acts on $E$ as translations.
Therefore any element of $G$ is of order $d=\deg L_{E}=1$.
So $X\cong D\times E\cong C\times F$, and $f:X\to C$ 
is the natural projection by construction.\qed
\enddemo
\proclaim{Lemma 1.5} Let $X$ be a smooth algebraic surface, 
$C$ a smooth curve, $f:X\to C$ a surjective morphism with connected fibers, 
and $F$ a general fiber of $f$. 
Then $q(X)\leq{g(C)+g(F)}$.
Moreover if this equality holds and $g(F)\geq 2$, 
then $X\sim_{bir} F\times C$.\endproclaim
\demo{Proof} 
See e.g. \cite{Be} p.345 or \cite{X}.\qed
\enddemo
\proclaim{Lemma 1.6}
Let $X$ be a minimal smooth surface of general type.
Then $K_{X}^{2}\geq 6q(X)-13$ unless $X\cong C_{1}\times C_{2}$
for some smooth curves $C_{1}$ and $C_{2}$.
\endproclaim
\demo{Proof}
We assume that $X\not\cong C_{1}\times C_{2}$ 
for smooth curves $C_{1}$ and $C_{2}$.
By Th\'eor\`eme 6.3 in \cite{D}, we have $K_{X}^{2}\geq 2p_{g}(X)+2(q(X)-4)+1$.
On the other hand, $p_{g}(X)\geq 2q(X)-3$ by \cite{Be}.
Hence $K_{X}^{2}\geq 6q(X)-13$. \qed
\enddemo
\proclaim{Proposition 1.7}
Let $X$ be a minimal smooth surface of general type 
and let $C$ be an irreducible reduced curve with $C^{2}>0$.
Then $K_{X}C\geq (3/2)q(X)-3$.
\endproclaim
\demo{Proof}If $q(X)\leq 2$, then this inequality is true.
So we assume $q(X)\geq 3$.
\flushpar
If $X\cong C_{1}\times C_{2}$ for some smooth curves $C_{1}$ and $C_{2}$,
then $K_{X}C\geq 2q(X)-4>(3/2)q(X)-3$.
So we may assume $X\not\cong C_{1}\times C_{2}$.
Let $x\in \Bbb{Q}$ with $x\geq 1$.
We put $m_{x}=g(xC)-q(X)$.
\proclaim{Claim 1.7.1}
If $2m_{x}+2\geq (2/3)(q(X)-2)+1$, then $(xC)^{2}\leq 2m_{x}+2$.
\endproclaim
\demo{Proof}
Assume that $(xC)^{2}>2m_{x}+2$.
Then $(xC)^{2}>(2/3)(q(X)-2)+1$.
\flushpar
Hence
$$\split
(K_{X})^{2}(xC)^{2}&>(6(q(X)-2)-1)(\frac{2}{3}(q(X)-2)+1) \\
                   &=4(q(X)-2)^{2}+6(q(X)-2)-\frac{2}{3}(q(X)-2)-1 \\
                   &=4(q(X)-2)^{2}+\frac{16}{3}(q(X)-2)-1
\endsplit
$$
by Lemma 1.6.
\flushpar
By Hodge index Theorem
$(xCK_{X})^{2}\geq (xC)^{2}(K_{X})^{2}>4(q(X)-2)^{2}$
and we have $xCK_{X}>2(q(X)-2)$.
Therefore 
$$
\split
g(xC)&>1+\frac{1}{2}(2(q(X)-2)+2m_{x}+2) \\
     &=q(X)+m_{x}
\endsplit
$$
and this is a contradiction.
\flushpar
This completes the proof of Claim 1.7.1.
\enddemo
We continue the proof of Proposition 1.7.
\flushpar
We have 
$$\split
q(X)+m_{x}=g(xC)&=g(C)+(x-1)g(C)+\frac{x-1}{2}(xC^{2}-2) \\
                &\geq q(X)+(x-1)q(X)+\frac{x-1}{2}(xC^{2}-2)
\endsplit
$$
since $g(C)\geq q(X)$.
\flushpar
Hence $m_{x}\geq (x-1)q(X)+((x-1)/2)(xC^{2}-2)$.
Here we put $x=(4/3)$.
Then $m_{x}\geq (1/3)q(X)-(1/9)>(1/3)q(X)-(7/6)$.
Therefore by Claim 1.7.1, we have 
$$(\frac{4}{3}C)^{2}\leq 2m_{x}+2.$$
In particular, $(4/3)CK_{X}\geq 2q(X)-4$.
Therefore $K_{X}C\geq (3/2)q(X)-3$.
This completes the proof of Proposition 1.7. \qed
\enddemo
\proclaim{Lemma 1.8}
Let $X$ be a minimal smooth surface of general type.
Then there are only finitely many irreducible curves $C$ on $X$ up to
numerical equivalence such that $K_{X}C$ is bounded.
\flushpar
Moreover there are only finitely many irreducible curves $C$ on $X$ 
such that $K_{X}C$ is bounded and $C^{2}<0$.
\endproclaim
\demo{Proof}
See Proposition 3 in \cite{Bo}.
\enddemo
\definition{Definition 1.9}(See e.g. \cite{BaBe}, \cite{BeFS}, and
\cite{BeS})
Let $X$ be a projective variety over $\Bbb{C}$ and
let $Z$ be a  0-dimensional subscheme of $X$.
A 0-dimensional subscheme $Z_{1}$ of $X$ is called a subcycle of $Z$ if 
$\Cal{I}_{Z}\subset\Cal{I}_{Z_{1}}$,
where $\Cal{I}_{Z}$ (resp. $\Cal{I}_{Z_{1}}$) is the ideal sheaf
which defines $Z$ (resp. $Z_{1}$).
Let $L$ be a Cartier divisor on $X$.  
Let $W$ be a subspace of $H^{0}(L)$ and $k$ a non-negative integer.
Then $W$ is called $k$-very ample
if the restriction map $W\to H^{0}(L\otimes\Cal{O}_{Z})$
is surjective for any 0-dimensional subscheme $Z$ with length $\leq k+1$.
If $W=H^{0}(L)$, then $LO$ is said to be k-very ample.
(We remark that $L$ is 0-very ample if and only if $L$ is spanned
and $L$ is 1-very ample if and only if $L$ is very ample.)
\enddefinition

\comment
\proclaim{Lemma 1.10} (Lemma 1.3 in \cite{Ba} or Lemma 2.8 in \cite{BaS})
Let $L$ be a line bundle on a smooth curve $C$.
Assume that a proper linear subspace $V\subset H^{0}(L)$ is $k$-very ample.
Then $\dim V\geq 2k+2$.
\endproclaim
\demo{Proof} See Lemma 1.3 in \cite{Ba}. \qed
\enddemo
\endcomment

\subhead{$\S$ 2. The case in which $\bold{\kappa(X)=0}$ or 1}
\endsubhead
\par
In this section, we will prove Conjecture 1$'$ for the case in
which $\kappa(X)=0$ or 1.
\proclaim{Theorem 2.1}
Let $(X,L)$ be a quasi-polarized surface with $\kappa(X)=0$ or 1.
Then $K_{X}L\geq 2q(X)-4$.
\flushpar
If this equality holds and $(X,L)$ is $L$-minimal, then
$(X,L)$ is one of the following;
\roster
\item $\kappa(X)=0$ case.\flushpar
$X$ is an Abelian surface and $L$ is any nef and big divisor.
\item $\kappa(X)=1$ case. \flushpar
$X\cong F\times C$ and $L\equiv C+(m+1)F$,
where $F$ and $C$ are smooth curves with $g(C)\geq 2$ and $g(F)=1$,
and $m=g(L)-q(X)$.
\endroster
\endproclaim
\demo{Proof}
(1) The case in which $\kappa(X)=0$.
\flushpar
Then $q(X)\leq 2$ by the classification theory of surfaces.
Hence $K_{X}L\geq 0\geq 2q(X)-4$.\flushpar
If $K_{X}L=2q(X)-4$, then $q(X)=2$ and $K_{X}L=0$.
Since $(X,L)$ is $L$-minimal,
$X$ is minimal, in particular, $X$ is an Abelian surface.
Conversely, let $(X,L)$ be any quasi-polarized surface which is $L$-minimal, 
and let $X$ be an Abelian surface.
Then $K_{X}L=0=2q(X)-4$.
\flushpar
(2) The case in which $\kappa(X)=1$.
\flushpar
Let $f:X\to C$ be an elliptic fibration,
$\mu :X\to X'$ the relatively minimal model of $X$,
and let $f': X'\to C$ be the relatively minimal elliptic fibration
such that $f=f'\circ \mu$.
Let $L'=\mu_{*}L$.
Then $L'$ is nef and big, and $K_{X}L\geq K_{X'}L'$.
\flushpar
By the canonical bundle formula for elliptic fibrations, we have
$$K_{X'}\equiv (2g(C)-2+\chi(\Cal{O}_{X'}))F'+\sum_{i}(m_{i}-1)F_{i},$$
where $F'$ is a general fiber of $f'$ 
and $m_{i}F_{i}$ is a multiple fiber of $f'$ for any $i$.
\flushpar 
Hence 
$$\split
K_{X'}L'\geq (2g(C)-2+\chi(\Cal{O}_{X'}))&\geq 2g(C)-2 \\
        &=2(g(C)+1)-4 \\
        &\geq 2q(X)-4.
\endsplit
$$
Therefore $K_{X}L\geq K_{X'}L'\geq 2q(X)-4$.
\flushpar
Assume that $K_{X}L=2q(X)-4$.
\flushpar
If $g(C)=0$, then $q(X)\leq 1$.
So $K_{X}L=2q(X)-4\leq -2$.
This is impossible.
Hence $g(C)\geq 1$.
By the above argument, we obtain $K_{X}L=K_{X'}L'=2q(X)-4$.
Since $(X,L)$ is $L$-minimal, $X$ is minimal.
Because $K_{X}L=2q(X)-4$, we obtain the following.
\roster
\item"(2-1)" $f$ has no multiple fibers.
\item"(2-2)" $\chi(\Cal{O}_{X})=0$.
\item"(2-3)" $q(X)=g(C)+1$.
\item"(2-4)" $LF=1$.
\endroster
By (2-3), (2-4), and Lemma 1.4,
we obtain $X\cong F\times C$ and $f: X\to C$ is the natural projection.
Because of $\kappa(X)=1$, we have $g(C)\geq 2$.
Then $f^{*}\circ f_{*}(\Cal{O}(L))\to \Cal{O}(L-Z)$ is surjective,
where $Z$ is a section of $f$.
Let $L|_{F_{t}}\sim p_{t}$, where $F_{t}=f^{-1}(t)$ and $t\in C$.
Let $(y,t)$ be a point of $F\times C$ and $(y(p_{t}),t)$ 
the point 
$p_{t}\in F\times C$.
Then the morphism $h: F\times C\to F\times C$ ;
$h(y,t)=(y-y(p_{t}),t)$ is an isomorphism.
Hence $L=h^{*}(\{ 0\}\times C)+f^{*}D$.
Therefore $L=C+f^{*}D$ via $h$, where $D\in \operatorname{Pic}(C)$.
But $L^{2}=2m+2$ for $m=g(L)-q(X)$.
Hence $L\equiv C+(m+1)F$.
This completes the proof of Theorem 2.1.\qed
\enddemo

\subhead{$\S$ 3. The case in which $\bold{\kappa(X)=2}$ and 
$\bold{h^{0}(L)\geq 2}$}
\endsubhead
\proclaim{Theorem 3.1}
Let $(X,L)$ be a quasi-polarized surface with $\kappa(X)=2$ and
$h^{0}(L)\geq 2$.
Then $K_{X}L\geq 2q(X)-2$.
\flushpar
If this equality holds and $(X,L)$ is $L$-minimal,
then $(X,L)$ is the following;
\flushpar
$X\cong F\times C$ and $L\equiv C+2F$,
where $F$ and $C$ are smooth curves with $g(F)=2$ and $g(C)\geq 2$.
\endproclaim
\demo{Proof}
(A) The case in which $X$ is minimal.
\flushpar
Then we use Theorem 1.3.
\flushpar
(A-1) The case in which $g(L)\geq 2q(X)-1$.
\flushpar
Then $q(X)+m=g(L)\geq 2q(X)-1$.
So we obtain $m\geq q(X)-1$.
\flushpar
(A-1-1) $q(X)\geq 1$ case.
\flushpar
Then by Lemma 1.2, we obtain $K_{X}^{2}\geq 2p_{g}(X)\geq 2q(X)$.
If $L^{2}\geq 2m$, then 
$$\split
(K_{X}L)^{2}\geq K_{X}^{2}L^{2}&\geq (2q(X))(2m) \\
                               &\geq 4q(X)(q(X)-1).
\endsplit
$$
Hence $K_{X}L>2(q(X)-1)$.
But this is impossible because
$$
\split
q(X)+m=g(L)&>1+\frac{1}{2}(2q(X)-2+2m) \\
           &=q(X)+m.
\endsplit
$$ 
Therefore $L^{2}\leq 2m-1$, that is, $K_{X}L\geq 2q(X)-1$.
\flushpar
(A-1-2) $q(X)=0$ case.
\flushpar
Then $K_{X}L>0>2q(X)-2$.
\flushpar
(A-2) The case in which $g(L)<2q(X)-1$.
\flushpar
Then by Theorem 1.3, there is a fiber space $f: X\to C$ such that 
$C$ is a smooth curve with $g(C)\geq 2$, $LF=1$, and $L-aF$ is numerically 
equivalent to an effective divisor,
where $F$ is a general fiber of $f$ and $a\geq 2$.
So there exists a section $C'$ of $f$ such that $C'$ is an irreducible
component of $L$, and we obtain that $L-aF\equiv C'+D'$,
where $D'$ is an effective divisor such that $f(D')$ are points.
\flushpar
Since $f$ is relatively minimal,
the relative canonical divisor $K_{X/C}=K_{X}-f^{*}K_{C}$ is nef
by Arakelov's Theorem.
So we have $K_{X/C}L\geq 2K_{X/C}F$.
Hence
$$
\split
g(L)&=g(C)+\frac{1}{2}(K_{X/C}L)+\frac{1}{2}L^{2} \\
    &\geq g(C)+K_{X/C}F+\frac{1}{2}L^{2} \\
    &=g(C)+2g(F)-2+\frac{1}{2}L^{2} \\
    &=g(C)+g(F)+g(F)-2+\frac{1}{2}L^{2} \\
    &\geq q(X)+\frac{1}{2}L^{2}
\endsplit
$$
because $g(F)\geq 2$ and $g(C)+g(F)\geq q(X)$.
\flushpar
Since $q(X)+m=g(L)$, we obtain $L^{2}\leq 2m$.
Namely $K_{X}L\geq 2q(X)-2$.\flushpar
Next we assume $K_{X}L=2q(X)-2$.
\flushpar
Then $g(L)<2q(X)-1$ by the above argument.
Moreover the following are satisfied by the above argument of (A-2);
\roster
\item"(a)" $K_{X/C}C'=0$, $K_{X/C}D'=0$.
\item"(b)" $a=2$.
\item"(c)" $g(F)=2$.
\item"(d)" $q(X)=g(C)+g(F)$.
\endroster
Since $X$ is minimal, we obtain $X\cong F\times C$ by (d) and
Lemma 1.5.
Moreover $f :X \to C$ is the natural projection.
Since $D'$ is contained in fibers of $f$ and $K_{X/C}D'=0$,
we obtain $D'=0$.
Since $K_{X/C}\equiv (2g(F)-2)C$ and $K_{X/C}C'=0$,
we have $CC'=0$.
Hence $C'$ is a fiber of $F\times C\to F$.
Therefore $L\equiv C+2F$ by (b).
\flushpar
(B) The case in which $X$ is not minimal.
\flushpar
Let $X=X_{0}\to X_{1}\to\cdots\to X_{n-1}\to X_{n}$ 
be the minimalization of $X$,
where $\mu_{i} :X_{i}\to X_{i+1}$ is the blowing down of (-1)-curve $E_{i}$.
Let $L_{i}=(\mu_{i-1})_{*}(L_{i-1})$,
$L_{0}=L$, and $L_{i-1}=(\mu_{i-1})^{*}L_{i}-m_{i-1}E_{i-1}$,
where $m_{i-1}\geq 0$.
We remark that $h^{0}(L_{i})=h^{0}((\mu_{i-1})^{*}L_{i})\geq h^{0}(L_{i-1})$.
Then $L^{2}=(L_{n})^{2}-\sum\limits_{i=0}^{n-1}m_{i}^{2}$
and $K_{X}L=K_{X_{n}}L_{n}+\sum\limits_{i=0}^{n-1}m_{i}$.
By the case (A), we have $K_{X_{n}}L_{n}\geq 2q(X)-2$.
Hence $K_{X}L\geq 2q(X)-2+\sum\limits_{i=0}^{n-1}m_{i}\geq 2q(X)-2$.
\flushpar
Next we consider $(X,L)$ such that $K_{X}L=2q(X)-2$ and
$(X,L)$ is $L$-minimal,
Then $\sum\limits_{i=0}^{n-1}m_{i}=0$ since $K_{X}L=2q(X)-2$
and so we have $m_{i}=0$ for any $i$.
But then $X$ is minimal because $(X,L)$ is $L$-minimal.
This is a contradiction.
This completes the proof of Theorem 3.1.\qed
\enddemo

\subhead{$\S$ 4. The case in which $\bold{\kappa(X)=2}$ and 
$\bold{h^{0}(L)=1}$}
\endsubhead
\par
In this section, we consider the case in which $\kappa(X)=2$
and $h^{0}(L)=1$.
We put $m=g(L)-q(X)$.
\proclaim{Lemma 4.1}
If $g(L)\geq 2q(X)$, then $K_{X}L\geq 2q(X)-1$.
\endproclaim
\demo{Proof}
Then $q(X)+m=g(L)\geq 2q(X)$.
Hence $m\geq q(X)$.
Assume that $L^{2}\geq 2m$.
So we obtain $L^{2}\geq 2q(X)$.
Let $\mu :X\to X'$ be the minimalization of $X$ and $L'=\mu_{*}L$.
Then $K_{X}L\geq K_{X'}L'$ and $(L')^{2}\geq L^{2}\geq 2q(X)$.
Since $K_{X'}^{2}\geq 2q(X)$ by Lemma 1.2,
we have
$(K_{X'}L')^{2}\geq (K_{X'})^{2}(L')^{2}\geq(2q(X))^{2}$
by Hodge index Theorem.
So we obtain $K_{X'}L'\geq 2q(X)$.
But this is impossible because
$$q(X)+m=g(L)\geq 1+q(X)+m.$$
Hence $L^{2}<2m$, that is, $K_{X}L\geq 2q(X)-1$.
This completes the proof of Lemma 4.1.\qed
\enddemo
\proclaim{Lemma 4.2}
If for any minimal quasi-polarized surfaces $(X,L)$ with $\kappa(X)=2$
and $h^{0}(L)\geq 1$ we can prove that $K_{X}L\geq 2q(X)-4$,
then this inequality holds for any quasi-polarized surface $(Y,A)$ with
$\kappa(Y)=2$ and $h^{0}(A)\geq 1$.
\endproclaim
\demo{Proof}
It is easy.\qed
\enddemo
\flushpar
By Lemma 4.2, it is sufficient to prove Conjecture 1 
(or Conjecture 1$'$) under the following assumption (4-1);
\flushpar
(4-1) $X$ is minimal.
\flushpar
Here we consider Conjecture 1 (or Conjecture 1$'$) for the following divisors.
\definition{Definition 4.3}
Let $X$ be a smooth projective surface and let $D$ be an effective divisor
on $X$.
Then $D$ is called a CNNS-divisor if the following conditions hold:
\roster
\item $D$ is connected.
\item the intersection matrix $\Vert (C_{i},C_{j})\Vert_{i,j}$
of $D=\sum\limits_{i}r_{i}C_{i}$ is not negative semidefinite.
\endroster
\enddefinition

\remark{Remark {\rm 4.4}}
If $L$ is an effective nef and big divisor, then $L$ is a CNNS-divisor.
\endremark

Let $D$ be a CNNS-divisor on a minimal smooth projective surface $X$ 
with $\kappa(X)=2$, 
and $D=\sum\limits_{i}r_{i}C_{i}$ its prime decomposition.
\flushpar
We divide three cases:
\roster
\item"($\alpha$)" $\sum\limits_{i\in S}r_{i}\geq 2$;
\item"($\beta$)"  $\sum\limits_{i\in S}r_{i}=1$;
\item"($\gamma$)" $\sum\limits_{i\in S}r_{i}=0$,
\endroster
where $S=\{\ i\ |\  C_{i}^{2}>0\}$.
\flushpar
First we consider the case ($\alpha$).
\proclaim{Theorem 4.5}
Let $D$ be a CNNS-divisor on a minimal smooth surface $X$ with $\kappa(X)=2$,
and let $D=\sum\limits_{i}r_{i}C_{i}$ be its prime decomposition.
If $\sum\limits_{i\in S}r_{i}\geq 2$, then $K_{X}L\geq 2q(X)-1$.
\endproclaim
\demo{Proof}
Let $A=\sum\limits_{i\in S}r_{i}C_{i}$ and $B=D-A$.
Then $A$ is nef and big.
We remark that $K_{X}D\geq K_{X}A$ since $X$ is minimal with $\kappa(X)=2$.
So it is sufficient to prove that $g(A)\geq 2q(X)$ by Lemma 4.1.
By assumption here, there are curves $C_{1}$ and $C_{2}$ 
(possibly $C_{1}=C_{2}$)
such that $C_{1}^{2}>0$ and $C_{2}^{2}>0$ and $A-C_{1}-C_{2}$ is effective.
Let $A_{12}=A-C_{1}-C_{2}$.
Then 
$$g(A)=g(C_{1}+C_{2})+\frac{1}{2}(K_{X}+A+C_{1}+C_{2})A_{12}.$$
Since $K_{X}+A$ is nef and $A$ is 1-connected,
we have $(K_{X}+A)A_{12}\geq 0$ and $(C_{1}+C_{2})A_{12}\geq 0$.
Hence $g(A)\geq g(C_{1}+C_{2})$.
On the other hand, 
$g(C_{1}+C_{2})=g(C_{1})+g(C_{2})+C_{1}C_{2}-1$.
Because $C_{1}^{2}>0$ and $C_{2}^{2}>0$,
we obtain $C_{1}C_{2}>0$.
Hence $g(C_{1}+C_{2})\geq g(C_{1})+g(C_{2})\geq 2q(X)$.
Therefore by Lemma 4.1,
we obtain $K_{X}(C_{1}+C_{2})\geq 2q(X)-1$.
So we have $K_{X}D\geq K_{X}(C_{1}+C_{2})\geq 2q(X)-1$.
This completes the proof of Theorem 4.5.\qed
\enddemo
Next we consider the case ($\gamma$).
\proclaim{Theorem 4.6}
Let $D$ be a CNNS-divisor on a minimal smooth projective surface $X$ 
with $\kappa(X)=2$ and let $D=\sum\limits_{i}r_{i}C_{i}$
be its prime decomposition.
If $\sum\limits_{i\in S}r_{i}=0$ and there exists a curve 
$C_{i}$ such that $C_{i}^{2}=0$, then $K_{X}D\geq 2q(X)-4$.
\endproclaim
\demo{Proof}
Assume that $C_{1}^{2}=0$.
We may assume that $q(X)\geq 1$.
Since $D$ is a CNNS-divisor, $D$ has at least two irreducible
components.
Let $C_{2}$ be another irreducible component of $D$
such that $C_{1}\cap C_{2}\not=\phi$.
Then 
$$g(D)=g(C_{1}+C_{2})+\frac{1}{2}(K_{X}+D+C_{1}+C_{2})D_{12},$$
where $D_{12}=D-(C_{1}+C_{2})$.
\flushpar
We put $l=g(C_{1}+C_{2})-q(X)$ and $m=g(D)-q(X)$.
Since $K_{X}D_{12}\geq 0$, we have $2m-2l\geq (D+C_{1}+C_{2})D_{12}$.
Let $X_{0}=X$, $C_{1,0}=C_{1}$, $C_{2,0}=C_{2}$,
and $\mu_{i} :X_{i}\to X_{i-1}$ blowing up at a point of 
$C_{1,i-1}\cap C_{2,i-1}$, where $C_{1,i}$ (resp. $C_{2,i}$) 
is the strict transform of $C_{1,i-1}$ (resp. $C_{2,i-1}$),
and let $E_{i}$ be an exceptional divisor such that $\mu_{i}(E_{i})$
is a point.
We put $\mu=\mu_{1}\circ\cdots\circ\mu_{n}$,
where $n$ is the natural number such that
$C_{1,n-1}\cap C_{2,n-1}\not=\phi$
and $C_{1,n}\cap C_{2,n}=\phi$.
Let $C_{1,i}=\mu_{i}^{*}C_{1,i-1}-b_{i}E_{i}$,
$C_{2,i}=\mu_{i}^{*}C_{2,i-1}-d_{i}E_{i}$,
and $a_{i}=b_{i}+d_{i}$.
We remark that $b_{i}\geq 1$ and $d_{i}\geq 1$.
Let $X_{0}'=X_{n}$, $C_{1,0}'=C_{1,n}$,
$C_{2,0}'=C_{2,n}$, $E_{0,0}'=E_{n}$,
and $\mu_{i}': X_{i}'\to X_{i-1}'$ blowing up
at a point 
$x\in (\operatorname{Sing}(C_{1,i-1}')\cup \operatorname{Sing}(C_{2,i-1}'))
\setminus ((C_{1,i-1}'\cap E_{0,i-1}')\cup (C_{2,i-1}'\cap E_{0,i-1}'))$,
where $C_{1,i}'$ (resp. $C_{2,i}'$, $E_{0,i}'$) is the
strict transform of $C_{1,i-1}'$ (resp. $C_{2,i-1}'$, $E_{0,i-1}'$),
and let $E_{i}'$ be an exceptional divisor on $X_{i}'$ such
that $\mu_{i}'(E_{i}')$ is a point.
Let $C_{1,i}'+C_{2,i}'=(\mu_{i}')^{*}(C_{1,i-1}'+C_{2,i-1}')-a_{i}'E_{i}'$.
We assume that
$(\operatorname{Sing}(C_{1,t}')\cup \operatorname{Sing}(C_{2,t}'))
\setminus ((C_{1,t}'\cap E_{0,t}')\cup (C_{2,t}'\cap E_{0,t}'))=\phi$.
\proclaim{Claim 4.7}
$g(C_{1,t}'+C_{2,t}'+E_{0,t}')\geq q(X_{t}')$.
\endproclaim
\demo{Proof}
Let $\alpha(C_{1,t}'+C_{2,t}'+E_{0,t}')=
\dim\operatorname{Ker}(H^{1}(\Cal{O}_{X_{t}'})\to 
H^{1}(\Cal{O}_{C_{1,t}'+C_{2,t}'+E_{0,t}'}))$.
By Lemma 3.1 in \cite{Fk4}, it is sufficient to prove 
$\alpha(C_{1,t}'+C_{2,t}'+E_{0,t}')=0$ since $C_{1,t}'+C_{2,t}'+E_{0,t}'$
is 1-connected.
Assume that $\alpha(C_{1,t}'+C_{2,t}'+E_{0,t}')\not=0$.
Since $q(X)\geq 1$, there is a morphism $f: X_{t}'\to G$
such that $f(X)$ is not a  point and $f(C_{1,t}'+C_{2,t}'+E_{0,t}')$ is a 
point,
where $G$ is an Abelian variety.
On the other hand, a $(\mu_{1}\circ\cdots\circ\mu_{n}\circ
\mu_{1}'\cdots\circ\mu_{t}')$-exceptional divisor is contracted by 
$f$ because $G$ is an Abelian variety.
Therefore $(\mu')^{*}(C_{1}+C_{2})$ is contracted by $f$.
But $(e_{1}C_{1}+C_{2})^{2}>0$ for sufficient large $e_{1}$.
This is impossible.
Hence $\alpha(C_{1,t}'+C_{2,t}'+E_{0,t}')=0$.
This completes the proof of Claim 4.7.
\enddemo
\flushpar
Hence 
$$
\split
g(C_{1,n}+C_{2,n}+E_{n})&=g(C_{1,t}'+C_{2,t}'+E_{0,t}')
                          +\sum_{i=1}^{t}\frac{1}{2}a_{i}'(a_{i}'-1) \\
                         &\geq q(X_{t}')
                          +\sum_{i=1}^{t}\frac{1}{2}a_{i}'(a_{i}'-1) \\
                         &=q(X)+\sum_{i=1}^{t}\frac{1}{2}a_{i}'(a_{i}'-1).
\endsplit
$$
On the other hand, 
$$g(C_{1}+C_{2})=g(C_{1,n}+C_{2,n}+E_{n})
+\sum_{i=1}^{n-1}\frac{1}{2}a_{i}(a_{i}-1)+\frac{1}{2}(a_{n}-1)(a_{n}-2).$$
Therefore
$$g(C_{1}+C_{2})\geq q(X)+\sum_{i=1}^{n-1}\frac{1}{2}a_{i}(a_{i}-1)
+\frac{1}{2}(a_{n}-1)(a_{n}-2)+\sum_{k=1}^{t}\frac{1}{2}a_{k}'(a_{k}'-1).$$
Since $l=g(C_{1}+C_{2})-q(X)$, we obtain
$$2l\geq \sum_{i=1}^{n-1}a_{i}(a_{i}-1)
+(a_{n}-1)(a_{n}-2)+\sum_{k=1}^{t}a_{k}'(a_{k}'-1).$$
Let $C_{1}C_{2}=x$.
Then $x=\sum\limits_{i=1}^{n}b_{i}d_{i}$
and $(C_{1}+C_{2})^{2}\leq 2x$ by hypothesis.
\proclaim{Claim 4.8}
$$2x-\sum_{i=1}^{n-1}a_{i}(a_{i}-1)-(a_{n}-1)(a_{n}-2)\leq 2.$$
\endproclaim
\demo{Proof}
$$
\split
&2x-\sum_{i=1}^{n-1}a_{i}(a_{i}-1)-(a_{n}-1)(a_{n}-2) \\
&=2\sum_{i=1}^{n}b_{i}d_{i}-\sum_{i=1}^{n-1}(b_{i}+d_{i})(b_{i}+d_{i}-1)-(b_{n}+d_{n}-1)(b_{n}+d_{n}-2).
\endsplit
$$
For each $i$($\not= n$),
$$\split
&2b_{i}d_{i}-(b_{i}+d_{i})(b_{i}+d_{i}-1) \\
&=-b_{i}^{2}-d_{i}^{2}+b_{i}+d_{i} \\
&=b_{i}(1-b_{i})+d_{i}(1-d_{i})\leq 0,
\endsplit
$$
and for $i=n$,
$$\split
&2b_{n}d_{n}-(b_{n}+d_{n}-1)(b_{n}+d_{n}-2) \\
&=-b_{n}^{2}-d_{n}^{2}+3b_{n}+3d_{n}-2 \\
&=b_{n}(3-b_{n})+d_{n}(3-d_{n})-2\leq 2.
\endsplit
$$
Therefore we obtain Claim 4.8.
\enddemo
By Claim 4.8, we obtain
$$
\split
D^{2}&=(C_{1}+C_{2})^{2}+(D+C_{1}+C_{2})D_{12} \\
     &\leq 2x+2m-2l \\
     &\leq 2x+2m-\sum_{i=1}^{n-1}a_{i}(a_{i}-1)
-(a_{n}-1)(a_{n}-2)-\sum_{k=1}^{t}a_{k}'(a_{k}'-1) \\
     &\leq 2m+2-\sum_{k=1}^{t}a_{k}'(a_{k}'-1) \\
     &\leq 2m+2.
\endsplit
$$
Therefore $K_{X}D\geq 2q(X)-4$.
This completes the proof of Theorem 4.6.\qed
\enddemo
Next we consider the case in which the equality in Theorem 4.6 holds.
\proclaim{Theorem 4.9}
Let $D$ be a CNNS-divisor on a minimal smooth surface $X$ with $\kappa(X)=2$,
and let $D=\sum\limits_{i}r_{i}C_{i}$ be its prime decomposition.
If $\sum\limits_{i\in S}r_{i}=0$, there exists a curve $C_{i}$
such that $C_{i}^{2}=0$, and $K_{X}D=2q(X)-4$,
then there are two irreducible curves $C_{1}$ and $C_{2}$
such that $D=C_{1}+C_{2}$ with $C_{1}^{2}=C_{2}^{2}=0$.
\flushpar
If $C_{1}$ or $C_{2}$ is not smooth, then $g(L)-q(X)=1$ or 3,
and $\sharp C_{1}\cap C_{2}=1$.
\roster
\item
If $g(L)-q(X)=1$, then $C_{i}$ is smooth but $C_{j}$ is 
not smooth only at $x\in C_{1}\cap C_{2}$
and $\operatorname{mult}_{x}(C_{j})=2$ for $i\not= j$
and $\{i,j\}=\{1,2\}$,
where $\operatorname{mult}_{x}(C_{j})$ is the multiplicity
of $C_{j}$ at $x$.
\item
If $g(L)-q(X)=3$, then $C_{1}$ and $C_{2}$ are  
not smooth only at $x\in C_{1}\cap C_{2}$
and $\operatorname{mult}_{x}(C_{i})=2$ for $i=$1, 2.
\endroster
\endproclaim
\demo{Proof}
Let $D=C_{1}+C_{2}+D_{12}$, where $C_{1}^{2}=0$ and $C_{2}$ is an
irreducible curve such that $C_{1}C_{2}>0$.
By the proof of Theorem 4.6, we have $K_{X}D_{12}=0$.
If $D_{12}\not=0$, then $K_{X}C=0$ for any irreducible curve $C$ of $D_{12}$
because $K_{X}$ is nef.
\proclaim{Claim 4.10}
$C^{2}=0$ for any irreducible curve $C$ of $D$.
\endproclaim
\demo{Proof}
By hypothesis, there is an irreducible curve $B$ of $D$
such that $B^{2}=0$.
Let $B'$ be any irreducible curve of $D$
such that $B\not= B'$ and $BB'>0$.
By the proof of Theorem 4.6 and $K_{X}D=2q(X)-4$,
we have $(B')^{2}=0$.
By repeating this argument, this completes the proof 
because $D$ is connected.
\enddemo
By this Claim, $C^{2}=0$ for any irreducible curve $C$ of $D_{12}$.
So $C\equiv 0$ by Hodge index Theorem.
But this is a contradiction.
\flushpar
Therefore $D_{12}=0$ and so we have $D=C_{1}+C_{2}$
with $C_{1}^{2}=C_{2}^{2}=0$.
Next we consider the singularity of $C_{1}$ and $C_{2}$.
\flushpar
We remark that $C_{1}$ (resp. $C_{2}$) is smooth
on $C_{1}\setminus \{C_{1}\cap C_{2}\}$
(resp. $C_{2}\setminus \{C_{1}\cap C_{2}\}$)
since $K_{X}D=2q(X)-4$ and $\sum\limits_{k=1}^{t}a_{k}'(a_{k}'-1)=0$
(here we use the notation in Theorem 4.6).
\flushpar
We assume that $\sharp C_{1}\cap C_{2}\geq 2$.
Then the number $n$ of blowing up $\mu=\mu_{1}\circ\cdots\circ\mu_{n}$ 
is greater than 1.
Since $K_{X}D=2q(X)-4$, we obtain $b_{1}=d_{1}=1$.
By interchanging the point of the first blowing up,
we obtain that $C_{1}$ and $C_{2}$ are smooth on $C_{1}\cap C_{2}$.
\flushpar
We assume $\sharp C_{1}\cap C_{2}=1$.
If the number $n$ of blowing up $\mu$ is greater than 1,
then $b_{1}=d_{1}=1$ by the proof of Theorem 4.6.
So $C_{1}$ and $C_{2}$ are smooth at $x\in C_{1}\cap C_{2}$.
Hence we assume that the number of blowing up is one.
Then $C_{1}C_{2}=b_{1}d_{1}$.
By the proof of Theorem 4.6,
$b_{1}(3-b_{1})+d_{1}(3-d_{1})=4$.
Hence $(b_{1},d_{1})=$(1,1), (1,2), (2,1), or (2,2).
\flushpar
If $(b_{1},d_{1})=(1,1)$, then $C_{1}$ and $C_{2}$ are smooth
at $x$.
\flushpar
If $(b_{1},d_{1})=$(1,2) or (2,1), then $C_{i}$ is smooth at $x$
and $C_{j}$ is not smooth at $x$ for $i\not=j$ and $\{i,j\}=\{1,2\}$,
and $\operatorname{mult}_{x}(C_{j})=2$,
where $\operatorname{mult}_{x}(C_{j})$ is the multiplicity
of $C_{j}$ at $x$.
In this case, $C_{1}C_{2}=2$ and $g(L)-q(X)=1$.
\flushpar
If $(b_{1},d_{1})=$(2,2), then $C_{1}$ and $C_{2}$ are 
not smooth at $x$,
and $\operatorname{mult}_{x}(C_{i})=2$ for $i=$1, 2.
In this case, $C_{1}C_{2}=4$ and $g(L)-q(X)=3$.
This completes the proof of Theorem 4.9.\qed
\enddemo
Next we consider the following case ($\ast$);
\roster
\item"($\ast$)"
Let $D$ be a CNNS-divisor on a minimal surface of general type,
and let $D=\sum\limits_{i}r_{i}C_{i}$ be its prime decomposition.
Then we assume $C_{i}^{2}<0$ for any $i$.
\endroster
\proclaim{Theorem 4.11}
Let (X,D) be ($\ast$).
Then $K_{X}D\geq 2q(X)-3$.
\endproclaim
Before we prove this theorem, we state some definitions and notaions
which is used in the proof of Theorem 4.11.
\definition{Definition 4.12}
Let $D$ be an effective divisor on $X$.
Then the dual graph $G(D)$ of $D$ is defined as follows.
\roster
\item The vertices of $G(D)$ corresponds to irreducible components of $D$.
\item For any two vertices $v_{1}$ and $v_{2}$ of $G(D)$, 
the number of edges joining $v_{1}$ and $v_{2}$ equal $\sharp\{ B_{1}\cap B_{2}\}$,
where $B_{i}$ is the component of $D$ corresponding to $v_{i}$ for $i=1$, 2.
\endroster
\enddefinition
\remark{Remark {\rm 4.12.1}}
Let $G(D)$ be the dual graph of an effective divisor $D$.
We reject one edge $e$ of $G(D)$ and $G=G(D)-\{e\}$.
Let $v_{1}$ and $v_{2}$ be vertices of $G(D)$ which are 
terminal points of the edge $e$.
Let $C_{1}$ and $C_{2}$ be the irreducible curve of $D$ 
corresponding $v_{1}$ and $v_{2}$ respectively.
Then $G$ is the dual graph
of the effective divisor which is the strict transform of $D$ 
by the blowing up at a point $x$ corresponding to $e$
if $i(C_{1},C_{2};x)=1$,
where $i(C_{i},C_{j};x)$ is the intersection number of $C_{i}$ and $C_{j}$
at $x$.
\endremark
\definition{Notation 4.13}
Let $(X,D)$ be ($\ast$).
We take a birational morphism $\mu': X'\to X$ such that
$C_{i}'\cap C_{j}'\cap C_{k}'=\phi$
for any distinct $C_{i}'$, $C_{j}'$, and $C_{k}'$,
and if $C_{i}'\cap C_{j}'\not=\phi$,
then $i(C_{i}',C_{j}';x)=1$ for $x\in C_{i}'\cap C_{j}'$,
where $D'=(\mu')^{*}(D)=\sum\limits_{i}r_{i}'C_{i}'$.
Let $\mu_{i}: X_{i}\to X_{i-1}$ be one point blowing up
such that $\mu'=\mu_{1}\circ\cdots\circ\mu_{t}$,
$X_{0}=X$ and $X_{t}=X'$.
Let $D_{i}=\mu_{i}^{*}D_{i-1}$ and $D_{0}=D$.
Let $b_{i}$ be an integer such that 
$(\mu_{i})^{*}((D_{i-1})_{\text{red}})-b_{i}E_{i}=(D_{i})_{\text{red}}$,
where $E_{i}$ is a $\mu_{i}$-exceptional curve.
\enddefinition
\remark{Remark {\rm 4.14}}
\flushpar
(a) No two $(\mu_{1}\circ\cdots\circ\mu_{i})$-exceptional curves on $X_{i}$
which are not (-1) curve intersect at a point on (-1)-curve on $X_{i}$
contracted by some $\mu_{j}$ ($j\leq i$).
\flushpar
(b) The point $x$ which is a center of blowing up $\mu_{i} :X_{i}\to X_{i-1}$
is contained in one of the following types;
\roster
\item the strict transform of the irreducible components of $D$;
\item the intersection of the strict transform of the irreducible 
components of $D$ and 
one (-1)-curve on $X_{i}$ contracted by some $\mu_{j}$ ($j\leq i$); 
\item
the intersection of the strict transform of the irreducible 
components of $D$ and 
one $(\mu_{1}\circ\cdots\circ\mu_{i})$-exceptional curve on $X_{i}$
which is not (-1)-curve and one (-1)-curve on $X_{i}$
contracted by some $\mu_{j}$ ($j\leq i$). 
\endroster
\endremark
We assume that $(X,D)$ satisfies ($\ast$) and we use Notation 4.13
unless specifically stated otherwise.
\definition{Definition 4.15}
\roster
\item Let $\pi : \widetilde{X}\to X$ be a birational morphism,
and let $\widetilde{X}$ and $X$ be smooth surfaces.
Let $\pi=\pi_{1}\circ\cdots \circ \pi_{n}$,
$X_{0}=X$, and $X_{n}=\widetilde{X}$, where $\pi_{i}:X_{i}\to X_{i-1}$
is a one point blowing up.
Let $E_{i}$ be the exceptional divisor of $\pi_{i}$.
Let $D$ be an effective divisor on $X$ and we put $D_{0}=D$.
Let $D_{i}=\pi^{*}(D_{i-1})$.
Then the multiplicity of $E_{i}$ in $D_{i}$ is called 
the $E_{i}$-multiplicity of $D$.
\item We use Notation 4.13.
Let $E$ be a union of $\mu$-exceptional curve and $D$ 
an effective divisor on $X$.
Let $x_{i}=\mu_{i}(E_{i})$.
If $x_{i}$ is the type (3) in Remark 4.14 (b), 
then a ($\mu_{1}\circ\cdots\circ\mu_{i}$)-exceptional curve $E$ 
which is not $(-1)$-curve is said to be an $e$-curve 
and $x_{i}$ is said to be an $e$-point.
\endroster
We remark that there is at most one $e$-curve throught $x_{i}$.
\enddefinition
\remark{Remark {\rm 4.16}}
We consider Notation 4.13.
Let $E$ an $e$-curve on $X_{i}$ and let $x_{i}$ be the $e$-point
associated with $E$.
Then we must be blowing up at $x_{i}$ by considering
Notation 4.13.
Let $\widetilde{E}$ be a strict transform of $E$ by
blowing up $\mu_{i+1}: X_{i+1}\to X_{i}$ at $x_{i}$.
Then $(\widetilde{E})^{2}=E^{2}-1\leq -3$
and $K_{X_{i+1}}\widetilde{E}=K_{X_{i}}E+1\geq 1$.
\endremark
\definition{Definition 4.17}
Let $\delta :\widetilde{X}\to X$ be any birational morphism,
$\widetilde{E}$ a union of $\delta$-exceptional curve,
and let $D$ be an effective divisor on $X$.
We put $B=\delta(\widetilde{E})=\{y_{1},\cdots,y_{s}\}$.
Then we can describe $\delta$ as $\delta=\delta_{s}\circ\cdots\circ\delta_{1}$,
where $\delta_{i}$ is the map whose image of a union of 
$\delta_{i}$-exceptional curves is $y_{i}$.
For each $y_{k}\in B$, we define a new graph $G=G(y_{k},D)$
which is called the {\it river} of the birational map $\delta_{k}$ 
and $D$.
\flushpar
(Step 1)
\flushpar
Let $E_{0,0}$ be a (-1)-curve obtained by blowing up at $y_{k}$.
Let $v_{0,0}$ be a vertex of the graph $G$ which corresponds to $E_{0,0}$.
We define the weight $u(0,0;G)$ of $v_{0,0}$ as follows:
$$u(0,0;G)=\text{the $E_{0,0}$-multiplicity of $D$.}$$
\flushpar
(Step 2)
\flushpar
Let $E_{1,1},\cdots, E_{1,t}$ be (-1)-curves
obtained by blowing up at distinct points
$\{x_{1,1},\cdots,x_{1,t}\}$ on $E_{0,0}$.
Let $v_{1,1},\cdots, v_{1,t}$ be vertices of the graph $G$ 
which correspond to $E_{1,1},\cdots, E_{1,t}$
respectively.
We join $v_{1,j}$ and $v_{0,0}$ by 
directed line which goes from $v_{1,j}$ to $v_{0,0}$.
For $j=1,\cdots, t$, we define the weight $u(1,j;G)$ of $v_{1,j}$
as follows:
$$u(1,j;G)=e_{1,j}-u(0,0;G),$$
where $e_{1,j}=\text{the $E_{1,j}$-multiplicity of $D$}$.
\flushpar
(Step 3)
\flushpar
In general, let $E_{i,1},\cdots, E_{i,t_{i}}$ be disjoint (-1)-curves
obtained by blowing up at distinct points
$\{x_{i,1},\cdots,x_{i,t_{i}}\}$ on $\bigcup\limits_{k}E_{i-1,k}$.
Let $v_{i,1},\cdots, v_{i,t_{i}}$ be vertices of the graph $G$ 
which correspond to $E_{i,1},\cdots, E_{i,t_{i}}$
respectively.
We join $v_{i,j}$ and $v_{i-1,k}$ by 
directed line which goes from $v_{i,j}$ to $v_{i-1,k}$
if $E_{i,j}$ 
is contracted in $E_{i-1,k}$.
Let $e_{i,j}=\text{the $E_{i,j}$-multiplicity of $D$}$
for $j=1,\cdots, t_{i}$.
Then we define the weight $u(i,j;G)$ of $v_{i,j}$
as follows:
$$u(i,j;G)=e_{i,j}-\sum_{v_{p,q}\in SP(i,j;G)}u(p,q;G),$$
where
$P(i,j;G)$ denotes the path between $v_{0,0}$ and $v_{i,j}$,
and $SP(i,j;G)=P(i,j;G)-\{v_{i,j}\}$.
\flushpar
By the above steps, we obtain the graph $G$ for each $y_{k}$.
\enddefinition
\definition{Notation 4.18}
$$w(i,j;G)=
\cases 
\deg(v_{i,j})-1, & \text{if $v_{i,j}\not=v_{0,0}$,} \\
\deg(v_{0,0}). & 
\endcases
$$
\enddefinition
\proclaim{Lemma 4.19}
Let $\mu :Y\to X$ be a birational morphism between
smooth surfaces $X$ and $Y$,
and let $D$ be an effective divisor on $X$.
Let $D'=\mu^{*}D$, and $E$ a union of all $\mu$-exceptional curves.
Let $B=\mu(E)$ and $M(D')=\text{ sum of the multiplicity of
(-1)-curve on $Y$ in $D'$}$.
Then
$$\split
M(D')&=\sum_{y\in B}\left[\sum_{v_{i,j}\in G(y)}\left
\{\sum_{v_{p,q}\in P(i,j;G(y))}u(p,q;G(y))\right\}\theta(i,j;G(y))\right]\\
&\ \ \ +\sum_{y\in B}\left\{\sum_{v_{i,j}\in G(y)}u(i,j;G(y))\right\},
\endsplit
$$
where $G(y)=G(y,D)$ and
$$\theta(i,j;G(y))=
\cases 
w(i,j;G(y))-1 & \text{if $w(i,j;G(y))\geq 1$,} \\
0 & \text{if $w(i,j;G(y))=0$.} 
\endcases
$$
\endproclaim
\demo{Proof}
We may assume that $B=\{y\}$.
Let $G=G(y,D)$.
Let $A=\{v_{i,j}\in G\ |\ \deg(v_{i,j})=1, v_{i,j}\not=v_{0,0}\}$
and $\rho=\sharp A-\deg(v_{0,0})$.
\flushpar
If $A=\phi$, then $M(D')=u(0,0;G)$.
\flushpar
So we asuume $A\not=\phi$.
We prove this Lemma by induction on the value of $\rho$.
We remark that by construction the following fact holds;
\proclaim{Fact}
For any $v_{s,t}\in A$,
the multiplicity of the (-1)-curve
corresponding to $v_{s,t}$ is equal to 
$\sum\limits_{v_{i,j}\in P(s,t;G)}u(i,j;G)$.
\endproclaim
\flushpar
(1) The case in which $\rho=0$.
\flushpar
Then $\deg{v}=2$ for any $v\not\in A$ and $v\not=v_{0,0}$.
Hence
$$\split
M(D')&=\sum_{v_{i,j}\in G}u(i,j;G)+u(0,0;G)(\deg(v_{0,0})-1) \\
     &=\sum_{v_{i,j}\in G}u(i,j;G)
    +\sum_{v_{i,j}\in G}\left\{ \sum_{v_{p,q}\in P(i,j;G)}u(p,q;G)\right\} 
     \theta(i,j;G).
\endsplit
$$
\flushpar
(2) The case in which $\rho=k>0$.
\flushpar
We assume that this Lemma is true for $\rho\leq k-1$.
We take a vertex $v_{s,t}\in A$ such that there is no edge 
whose terminal points are
$v_{0,0}$ and $v_{s,t}$.
Let $G^{\vee}=G-\{v_{s,t}\}$.
Let $\mu^{-} :Y\to X^{-}$ be blowing down of (-1)-curves corresponding to
$v_{s,t}$ and $\mu=\mu^{+}\circ\mu^{-}$.
Let $D^{\vee}=(\mu^{+})^{*}(D)$.
Then we remark that $G^{\vee}$ is the river of $\mu^{+}$ and $D$.
\flushpar
Then by induction hypothesis
$$M(D^{\vee})=\sum_{v_{i,j}\in G^{\vee}}u(i,j;G^{\vee})
+\sum_{v_{i,j}\in G^{\vee}}\left\{ \sum_{v_{p,q}\in P(i,j;G^{\vee})}
u(p,q;G^{\vee})\right\} \theta(i,j;G^{\vee}).
$$
Next we consider $M(D')$.
Let $v_{s-1,l}$ be a vertex such that there is an edge
between $v_{s-1,l}$ and $v_{s,t}$.
\flushpar
(2-1) The case in which $w(s-1,l;G)=1$.
\flushpar
Then $M(D')=M(D^{\vee})+u(s,t;G)$.
Hence 
$$
\split
M(D')=&\sum_{v_{i,j}\in G^{\vee}}u(i,j;G^{\vee})+u(s,t;G)
+\sum_{v_{i,j}\in G^{\vee}}\left\{ \sum_{v_{p,q}\in P(i,j;G^{\vee})}
u(p,q;G^{\vee})\right\} \theta(i,j;G^{\vee}) \\
&=\sum_{v_{i,j}\in G}u(i,j;G)
+\sum_{v_{i,j}\in G}\left\{\sum_{v_{p,q}\in P(i,j;G)}
u(p,q;G)\right\}\theta(i,j;G),
\endsplit
$$
\flushpar
because $\theta(s-1,l;G)=\theta(s,t;G)=0$
and we have
$u(i,j;G)=u(i,j;G^{\vee})$,
$w(i,j;G)=w(i,j;G^{\vee})$, and $\theta(i,j;G)=\theta(i,j;G^{\vee})$
for $v_{i,j}\not=v_{s,t}$ .
\flushpar
(2-2) The case in which $w(s-1,l;G)\geq 2$.
\flushpar
Then 
$$M(D')=M(D^{\vee})+\sum_{v_{p,q}\in SP(s,t;G)}u(p,q;G)+u(s,t;G).$$
Hence
$$\split
M(D')&=\sum_{v_{i,j}\in G^{\vee}}u(i,j;G^{\vee})+u(s,t;G)
+\sum_{v_{i,j}\in G^{\vee}}\left\{\sum_{v_{p,q}\in P(i,j;G^{\vee})}
u(p,q;G^{\vee})\right\} \theta(i,j;G^{\vee}) \\
&\ \ \ +\sum_{v_{p,q}\in SP(s,t;G)}u(p,q;G) \\
&=\sum_{v_{i,j}\in G}u(i,j;G)
+\sum_{v_{i,j}\in G}\left\{ \sum_{v_{p,q}\in P(i,j;G)}
u(p,q;G)\right\} \theta(i,j;G),
\endsplit
$$
\flushpar
because $\theta(s,t;G)=0$ and $\theta(s-1,l;G)=\theta(s-1,l;G^{\vee})+1$
and because we have 
$u(i,j;G)=u(i,j;G^{\vee})$,
$w(i,j;G)=w(i,j;G^{\vee})$, and $\theta(i,j;G)=\theta(i,j;G^{\vee})$
for $(i,j)\not=(s,t),(s-1,l)$.
This completes the proof of Lemma 4.19.\qed
\enddemo
\proclaim{Lemma 4.20}
Let $D$ be a CNNS-divisor on $X$ and we use Notation 4.13.
Then 
$$(D_{\text{red}}')^{2}\leq 2l-2-\sum_{i=1}^{t}b_{i}(b_{i}-1)
+\sum_{j}((C_{j}')^{2}+2).$$
\endproclaim
\demo{Proof}
First we prove the following Claim.
\proclaim{Claim 4.21}
$$e(D')-o(D')+1+\sum_{i=1}^{t}\frac{1}{2}b_{i}(b_{i}-1)\leq l,$$
where $l=g(D_{\text{red}})-q(X)$.
\endproclaim
\demo{Proof}
We have 
$g(D_{\text{red}}')=
g(D_{\text{red}})-\sum\limits_{i=1}^{t}\frac{1}{2}b_{i}(b_{i}-1)$ 
by definition.
There exists $m=e(D')-o(D')+1$ egdges $e_{1},\cdots,e_{m}$
of $G(D_{\text{red}})$ such that $G-\{e_{1},\cdots,e_{m}\}$
is a tree.
Therefore by Remark 4.12.1, there exists a connected effective divisor 
$A$ on $X''$ which is obtained by finite number of blowing ups 
of $X'$ such that $g(D_{\text{red}}')=g(A)+e(D')-o(D')+1$.
Let $\mu'' :X''\to X'$ be its birational morphism
and $A$ the strict transform of $D_{\text{red}}'$ by $\mu''$.
Let $\alpha(A)=\dim\operatorname{Ker}(H^{1}(\Cal{O}_{X''})
\to H^{1}(\Cal{O}_{A}))$.
Then we calculate $\alpha(A)$.
\flushpar
If $\alpha(A)\not=0$,
then there exist an Abelian variety $T$, a surjective morphism
$f'':X''\to T$ such that $f''(X'')$ is not a point and $f''(A)$ is a point.
Then any $\mu''$-exceptional curve is contracted by $f''$ 
because $T$ is an Abelian variety.
Hence $f''((\mu'')^{*}D_{\text{red}}')$ is a point.
But $(\mu'')^{*}D_{\text{red}}'$ is not negative semidefinite.
Therefore $\alpha(A)=0$.
Since $A$ is reduced and connected, $A$ is 1-connected.
Hence $g(A)=h^{1}(\Cal{O}_{A})$.
So we obtain $g(A)=h^{1}(\Cal{O}_{A})\geq q(X'')=q(X)$.
\flushpar
By the above argument,
$$\split
g(D_{\text{red}})&=g(D_{\text{red}}')
                   +\sum_{i=1}^{t}\frac{1}{2}b_{i}(b_{i}-1) \\
                 &=g(A)+e(D')-o(D')+1
                   +\sum_{i=1}^{t}\frac{1}{2}b_{i}(b_{i}-1) \\
                 &\geq q(X)+e(D')-o(D')+1
                   +\sum_{i=1}^{t}\frac{1}{2}b_{i}(b_{i}-1). 
\endsplit
$$
Therefore
$$e(D')-o(D')+1+\sum_{i=1}^{t}\frac{1}{2}b_{i}(b_{i}-1)\leq l.$$               This completes the proof of Claim 4.21.
\enddemo
We continue the proof of Lemma 4.20.
By construction, we obtain
$$\split
(D_{\text{red}}')^{2}&=\sum_{j}(C_{j}')^{2}+2e(D') \\
                     &=\sum_{j}(C_{j}')^{2}+2(o(D')+e(D')-o(D')) \\
                     &=\sum_{j}((C_{j}')^{2}+2)+2(e(D')-o(D')).
\endsplit
$$
By Claim 4.21, we have
$$(D_{\text{red}}')^{2}\leq 2l-2-\sum_{i=1}^{t}b_{i}(b_{i}-1)
+\sum_{j}((C_{j}')^{2}+2).$$
This completes the proof of Lemma 4.20.\qed
\enddemo
\proclaim{Theorem 4.22}
Let $X$ be a minimal smooth projective surface with $\kappa(X)\geq 0$
and $D$ a CNNS-divisor on $X$.
Let $D=\sum\limits_{j}r_{j}D_{j}$ be its prime decomposition
and $m=g(D)-q(X)$,
where $m\in\Bbb{Z}$.
\flushpar
Then $D^{2}\leq 2m-2+N(D)$,
where 
$$N(D)=\sum\limits_{\beta\in\Bbb{Z}}\beta\cdot\sharp
\{ \text{irreducible curves $C_{j}$ of $D$ such that $C_{j}^{2}=-2+\beta$}\}.$$
\endproclaim
\demo{Proof}
We use Notation 4.13 and the notions which is defined above.
We may assume that $B=\{y\}$.
Let $G=G(y,D)$, $u(i,j)=u(i,j;G)$,
$\theta(i,j)=\theta(i,j;G)$,
$w(i,j)=w(i,j;G)$, $P(i,j)=P(i,j;G)$, 
and $SP(i,j)=SP(i,j;G)$.
Let $D'=(\mu')^{*}D$ and $D_{\text{nr}}'=D'-D_{\text{red}}'$.
Let $D_{\text{nr}}'=D_{\text{ne}}'+D_{\text{e}}'+D_{-1}'$,
where $D_{\text{ne}}'$ is the effective divisor which consists of not 
$\mu'$-exceptional curves,
$D_{\text{e}}'$ is the effective divisor which consists of curves
which are $\mu'$-exceptional curves and not (-1)-curves,
and $D_{-1}'$ is the effective divisor which consists of (-1)-curves.
\flushpar
Then
$$K_{X'}D_{\text{e}}'=\sum_{v_{i,j}\in G}\left\{ 
\left (\sum_{v_{p,q}\in P(i,j)}u(p,q)\right )-1\right\}\theta(i,j)
+\sum_{v_{i,j}\in G}\varepsilon(i,j)(m(i,j)-1),$$
where $m(i,j)$ is the multiplicity of $e$-curve through $x_{i,j}$ 
in the total transform of $D$, 
$x_{i,j}$ is the blowing up point and its (-1)-curve corresponds to $v_{i,j}$,
$\varepsilon(i,j)=1$ if there exists the $e$-curve through $x_{i,j}$
and $\varepsilon(i,j)=0$ if there does not exist the $e$-curve through
$x_{i,j}$.
\flushpar
On the other hand,
$$-\sum_{\alpha}(E_{\alpha}^{2}+2)=\sum_{v_{i,j}\in G-W}(w(i,j)-1)
+\sum_{v_{i,j}\in G}\varepsilon(i,j),$$
where $E_{\alpha}$ is a $\mu'$-exceptional curve on $X'$ and not (-1)-curve,
and $W=\{v_{i,j}\in G| w(i,j)=0\}$.
\flushpar
Hence
$$\align
K_{X'}D_{\text{e}}'-\sum_{\alpha}(E_{\alpha}^{2}+2)&=
\sum_{v_{i,j}\in G}\left\{ 
\left (\sum_{v_{p,q}\in P(i,j)}u(p,q)\right )-1\right\}\theta(i,j) 
\tag 4.22.1 \\
&\ \ \ +\sum_{v_{i,j}\in G}\varepsilon(i,j)m(i,j)
+\sum_{v_{i,j}\in G-W}(w(i,j)-1).
\endalign
$$
Let 
$$\beta_{\text{nr}}=\text{sum of multiplicity of $\mu'$-exceptional 
(-1)-curves in $D_{\text{nr}}'$}.$$
Then
$$\align -\beta_{\text{nr}}=K_{X'}D_{-1}'. \tag 4.22.2
\endalign
$$
Let $C_{i,j}$ be a strict transform of $C_{i,j-1}$
by $\mu_{j}$ and $C_{i,0}=C_{i}$.
Let $C_{i,j}=\mu_{j}^{*}(C_{i,j-1})-e(i)_{j}E_{j}$,
where $E_{j}$ is the (-1)-curve of $\mu_{j}$.
We remark that $e(i)_{j}\geq 1$ for any $i$, $j$.
\flushpar
Then
$$K_{X'}((r_{i}-1)C_{i,t})\geq \sum_{i=1}^{t}(r_{i}-1)e(i)_{j}$$
because $X$ is minimal.
\flushpar
Hence 
$$K_{X'}(D_{\text{ne}}')\geq \sum_{i}\left\{\sum_{j=1}^{t}(r_{i}-1)e(i)_{j}
\right\}.$$
On the other hand
$$\sum_{i}(C_{i,t}^{2}+2)=N(D)-\sum_{i}\sum_{j=1}^{t}e(i)_{j}^{2}$$
\flushpar
because $C_{i,t}^{2}=C_{i}^{2}-\sum\limits_{j=1}^{t}e(i)_{j}^{2}$.
\flushpar
Hence 
$$\align
K_{X'}(D_{\text{ne}}')-\sum_{i}(C_{i,t}^{2}+2) 
\geq \sum_{i}\sum_{j=1}^{t}(r_{i}e(i)_{j})-N(D) \tag 4.22.3 
\endalign
$$
\flushpar
since $\sum\limits_{j=1}^{t}e(i)_{j}^{2}\geq \sum\limits_{j=1}^{t}e(i)_{j}$.
\flushpar
By  (4.22.1), (4.22.2), and (4.22.3),
we obtain
$$\align
&K_{X'}D_{\text{nr}}'-\sum_{i}(C_{i,t}^{2}+2)-\sum_{\alpha}(E_{\alpha}^{2}+2)
\tag 4.22.4 \\
&\geq -\beta_{\text{nr}}+\sum_{v_{i,j}\in G}\left\{ 
\left (\sum_{v_{p,q}\in P(i,j)}u(p,q)\right )-1\right\}\theta(i,j) \\
&\ \ \ +\sum_{v_{i,j}\in G}\varepsilon(i,j)m(i,j)
+\sum_{v_{i,j}\in G-W}(w(i,j)-1)+\sum_{i}\sum_{j=1}^{t}(r_{i}e(i)_{j})-N(D).
\endalign
$$
On the other hand, we have
$$\split
q(X)+m=g(D)&=g(D') \\
&=g(D_{\text{red}}')+\frac{1}{2}(K_{X'}+D'+D_{\text{red}}')D_{\text{nr}}'\\
&=g(D_{\text{red}})-\frac{1}{2}\sum_{i=1}^{t}b_{i}(b_{i}-1)
+\frac{1}{2}(K_{X'}+D'+D_{\text{red}}')D_{\text{nr}}'\\
&=q(X)+l-\frac{1}{2}\sum_{i=1}^{t}b_{i}(b_{i}-1)
+\frac{1}{2}(K_{X'}+D'+D_{\text{red}}')D_{\text{nr}}',
\endsplit
$$
where $l=g(D_{\text{red}})-q(X)$.
\flushpar
Hence by (4.22.4), we obtain
$$\split
2m-2l&=(K_{X'}+D'+D_{\text{red}}')D_{\text{nr}}'-\sum_{i=1}^{t}b_{i}(b_{i}-1)\\
&\geq \sum_{i}(C_{i,t}^{2}+2)+\sum_{\alpha}(E_{\alpha}^{2}+2)
-\beta_{\text{nr}} \\
&\ \ \ +\sum_{v_{i,j}\in G}\left\{ 
\left (\sum_{v_{p,q}\in P(i,j)}u(p,q)\right )-1\right\}\theta(i,j) \\
&\ \ \ +\sum_{v_{i,j}\in G}\varepsilon(i,j)m(i,j)
+\sum_{v_{i,j}\in G-W}(w(i,j)-1) \\
&\ \ \ +\sum_{i}\sum_{j=1}^{t}(r_{i}e(i)_{j})-N(D)
+(D'+D_{\text{red}}')D_{\text{nr}'}-\sum_{i=1}^{t}b_{i}(b_{i}-1),
\endsplit
$$
\flushpar
and so we have
$$\split
(D'+D_{\text{red}}')D_{\text{nr}'}\leq 
& -\sum_{i}(C_{i,t}^{2}+2)-\sum_{\alpha}(E_{\alpha}^{2}+2)
+\beta_{\text{nr}}\\
&\ \ \ -\sum_{v_{i,j}\in G}\left\{ 
\left (\sum_{v_{p,q}\in P(i,j)}u(p,q)\right )-1\right\}\theta(i,j) \\
&\ \ \ -\sum_{v_{i,j}\in G}\varepsilon(i,j)m(i,j)
-\sum_{v_{i,j}\in G-W}(w(i,j)-1) \\
&\ \ \ -\sum_{i}\sum_{j=1}^{t}(r_{i}e(i)_{j})+N(D)
+\sum_{i=1}^{t}b_{i}(b_{i}-1)+2m-2l.
\endsplit
$$
Therefore by Lemma 4.20, we obtain
$$
\split
(D')^{2}&=(D_{\text{red}}')^{2}+(D'+D_{\text{red}}')D_{\text{nr}}' \\
        & \leq (2m-2l)+(2l-2)+\sum_{i=1}^{t}b_{i}(b_{i}-1)
         -\sum_{i=1}^{t}b_{i}(b_{i}-1) \\
        & \ \ \ +\sum_{i}((C_{i}')^{2}+2)-\sum_{i}(C_{i,t}^{2}+2)
         -\sum_{\alpha}(E_{\alpha}^{2}+2)+\beta_{\text{nr}} \\
        & \ \ \ -\sum_{v_{i,j}\in G}\left\{ 
          \left (\sum_{v_{p,q}\in P(i,j)}u(p,q)\right)-1\right\}\theta(i,j) \\
        & \ \ \ -\sum_{v_{i,j}\in G}\varepsilon(i,j)m(i,j)
           -\sum_{v_{i,j}\in G-W}(w(i,j)-1) \\
        & \ \ \  -\sum_{i}\sum_{j=1}^{t}(r_{i}e(i)_{j})+N(D) \\
        &=(2m-2)+M(D') \\
        & \ \ \ -\sum_{v_{i,j}\in G}\left\{ 
          \left (\sum_{v_{p,q}\in P(i,j)}u(p,q)\right )-1\right\}\theta(i,j) \\
        & \ \ \ -\sum_{v_{i,j}\in G}\varepsilon(i,j)m(i,j)
           -\sum_{v_{i,j}\in G-W}(w(i,j)-1) \\
        & \ \ \  -\sum_{i}\sum_{j=1}^{t}(r_{i}e(i)_{j})+N(D),
\endsplit
$$
where $M(D')$ is the sum of the multiplicity of (-1)-curves in $D'$.
\flushpar
On the other hand by Lemma 4.19, we have
$$
\split
&M(D') -\sum_{v_{i,j}\in G}\left\{ 
          \left (\sum_{v_{p,q}\in P(i,j)}u(p,q)\right)-1\right\}\theta(i,j)\\
&=M(D') -\sum_{v_{i,j}\in G}\left\{ 
          \sum_{v_{p,q}\in P(i,j)}u(p,q)\right\}\theta(i,j)
  +\sum_{v_{i,j}\in G-W}(w(i,j)-1) \\
&= \sum_{v_{p,q}\in G}u(p,q)+\sum_{v_{i,j}\in G-W}(w(i,j)-1).
\endsplit
$$
Therefore
$$\split
(D')^{2}&\leq 2m-2+\sum_{v_{p,q}\in G}u(p,q)+\sum_{v_{i,j}\in G-W}(w(i,j)-1) \\
& \ \ \ -\sum_{v_{i,j}\in G}\varepsilon(i,j)m(i,j)
           -\sum_{v_{i,j}\in G-W}(w(i,j)-1) \\
& \ \ \  -\sum_{i}\sum_{j=1}^{t}(r_{i}e(i)_{j})+N(D)\\
& =2m-2+N(D)
\endsplit
$$
because we have 
$$\sum_{v_{i,j}\in G}\varepsilon(i,j)m(i,j)
+\sum_{i}\sum_{j=1}^{t}(r_{i}e(i)_{j})=\sum_{v_{p,q}\in G}u(p,q)$$
\flushpar
by considering the definition of $u(p,q)$.
This completes the proof of Theorem 4.22. \qed
\enddemo
Theorem 4.11 is obtained by Theorem 4.22.
\demo{Proof of Theorem 4.11}
It is sufficient to prove $D^{2}\leq 2m+1$
if $g(D)-q(X)=m$.
We consider the following decomposition ($\ast\ast$) of $D$;
\roster
\item"($\ast\ast$)"
$D=D_{1}+D_{2}$, and $D_{1}$ and $D_{2}$ have no common
component, where $D_{1}$ and $D_{2}$ are non zero effective connected
divisors,
\endroster
\proclaim{Claim 4.23}
If $((D_{1})_{\text{red}})^{2}\leq 0$ and
$((D_{2})_{\text{red}})^{2}\leq 0$, then $N(D)\leq 4$.
\flushpar
If $((D_{1})_{\text{red}})^{2}<0$ or
$((D_{2})_{\text{red}})^{2}<0$, then $N(D)\leq 3$.
\endproclaim
\demo{Proof}
Let $(D_{i})_{\text{red}}=\sum\limits_{j}B_{i,j}$.
Then $\sum\limits_{j}(B_{i,j})^{2}=N(D_{i})-2o(D_{i})$
and $\sum\limits_{j\not=k}B_{i,j}B_{i,k}\geq e(D_{i})$.
Hence 
$((D_{i})_{\text{red}})^{2}\geq 2e(D_{i})-2o(D_{i})+N(D_{i})$
for $i=1$, 2.
By hypothesis, we have $0\geq 2e(D_{i})-2o(D_{i})+N(D_{i})$
for $i=$1, 2.
Since the dual graph $G(D_{i})$ of $D_{i}$ is connected, we have
$e(D_{i})-o(D_{i})+1\geq 0$.
Hence $2e(D_{i})-2o(D_{i})\geq -2$
and so we have $N(D_{i})\leq 2$.
\flushpar
On the other hand, $N(D)=N(D_{1})+N(D_{2})$
since $D=D_{1}+D_{2}$.
Therefore $N(D)\leq 4$.
\flushpar
The last part of Claim 4.23 can be proved by the above argument.
This completes the proof of Claim 4.23.
\enddemo
Let $S(D)$ be a set of an effective connected reduced divisor 
$\widetilde{D}$ contained in $D$ such that $\widetilde{D}$
has a minimum component which satisfies the property that
the intersection matrix of $\widetilde{D}$ is not negative semidefinite.
\flushpar
Then $S(D)\not=\phi$ by hypothesis.
Let $\bar{D}=\sum\limits_{i\in J}C_{i}\in S(D)$
and let $r_{i}$ be the multiplicity of $C_{i}$ in $D$.
Let $D_{\alpha}=\sum\limits_{i\in J}r_{i}C_{i}$
and $D_{\beta}=D-D_{\alpha}$.
We remark that possibly $D_{\beta}=0$.
Then $D_{\alpha}$ has at least two components since $C_{i}^{2}<0$ for any $i$.
Let $D_{\alpha}=D_{\alpha ,1}+D_{\alpha ,2}$
be the decomposition as ($\ast\ast$).
\proclaim{Claim  4.24}
We can take this decomposition which satisfies $(D_{\alpha ,1})^{2}<0$.
\endproclaim
\demo{Proof}
We consider the dual graph $G(D_{\alpha})$ of $D_{\alpha}$.
Then $G(D_{\alpha})$ is connected.
In Graph Theory, there is the following standard Theorem;
\proclaim{Theorem 4.25}
Let $G$ be a connected graph which is not one point.
Then there are at least two points which are not cutpoints.
(Here a vertex $v$ of a graph is called a cutpoint
if removal of $v$ increases the number of components.)
\endproclaim
\demo{Proof}
See Theorem 3.4 in \cite{H}.
\enddemo
We continue the proof of Claim 4.24.
By Theorem 4.25, it is sufficient to take $(D_{\alpha ,1})_{\text{red}}$
as an irreducible curve corresponding to a vertex of $G(D_{\alpha})$
which is not a cutpoint.
This completes the proof of Claim 4.24.
\enddemo
We continue the proof of Theorem 4,11.
\flushpar
We have $((D_{\alpha ,1})_{\text{red}})^{2}<0$
and $((D_{\alpha ,2})_{\text{red}})^{2}\leq 0$
by the choice of $D_{\alpha}$.
Therefore $N(D_{\alpha})\leq 3$ by Claim 4.23.
\flushpar
On the other hand, we have
$$q(X)+m=g(D)=g(D_{\alpha})+\frac{1}{2}(K_{X}+D+D_{\alpha})D_{\beta}.$$
Let $g(D_{\alpha})=q(X)+m_{\alpha}$.
Then by Theorem 4.22, $D_{\alpha}^{2}\leq 2m_{\alpha}-2+N(D_{\alpha})\leq 2m_{\alpha}+1$
since $D_{\alpha}$ is a CNNS-divisor.
\flushpar
On the other hand, $(K_{X}+D+D_{\alpha})D_{\beta}=2(m-m_{\alpha})$
and $K_{X}D_{\beta}\geq 0$.
Hence $(D+D_{\alpha})D_{\beta}\leq 2(m-m_{\alpha})$.
Therefore
$$\split
D^{2}&=D_{\alpha}^{2}+(D+D_{\alpha})D_{\beta} \\
     &\leq 2m_{\alpha}+1+2m-2m_{\alpha} \\
     &=2m+1.
\endsplit
$$
This completes the proof of Theorem 4.11. \qed
\enddemo
\remark{Remark {\rm 4.26}}
Let $D=\sum\limits_{i}r_{i}C_{i}$ be an effective divisor
on a minimal smooth surface of general type with $C_{i}^{2}<0$ for any $i$.
If the intersection matrix $\Vert (C_{i}\cdot C_{j})\Vert$
is not negative semidefinite, then $K_{X}D\geq 2q(X)-3$.
\flushpar
Indeed, let $D_{1},\cdots, D_{t}$ be the connected component of $D$.
Then for some $D_{k}$, the intersection matrix of the
components of $D$ is not negative semidefinite.
By Theorem 4.11, we have $K_{X}D_{k}\geq 2q(X)-3$.
Since $K_{X}$ is nef, we obtain $K_{X}D\geq 2q(X)-3$.
\endremark
\proclaim{Corollary 4.27}
Let $X$ be a smooth surface of general type
and let $D$ be a nef-big effective divisor with $h^{0}(D)=1$
on $X$.
If $D$ is not the following type ($\star$), then
$K_{X}D\geq 2q(X)-4$;
\roster
\item"($\star$)"
$D=C_{1}+\sum\limits_{j\geq 2}r_{j}C_{j}$ ;
$C_{1}^{2}>0$ and the intersection matrix 
$\Vert C_{j}, C_{k}\Vert_{j\geq 2, k\geq 2}$
of $\sum\limits_{j\geq 2}r_{j}C_{j}$ is
negative semidefinite.
\endroster
\endproclaim
\demo{Proof}
By Theorem 4.5, Theorem 4.6, Theorem 4.11, and Remark 4.26,
we obtain Corollary 4.27. \qed
\enddemo

\subhead{$\S$ 5. The case in which $\bold{\kappa(X)=2}$ and 
$\bold{L}$ is an irreducible 
reduced curve}
\endsubhead
\par
\definition{Notation 5.1}
Let $X$ be a smooth projective surface over the complex number field $\Bbb{C}$
and let $C$ be a curve on $X$ with $C^{2}>0$.
Let $N(k;C)$ be the set of a 0-dimensional subscheme $\widetilde{Z}$
with $\operatorname{length}\widetilde{Z}=k+1$ and
$\operatorname{Supp}\widetilde{Z}\subset C$ such that the restriction map
$\Gamma(\Cal{O}(K_{X}+C))\to 
\Gamma(\Cal{O}(K_{X}+C)\otimes\Cal{O}_{\widetilde{Z}})$
is not surjective.
Let $S(\widetilde{Z};C)$ be the set of a subcycle $Z$ of 
$\widetilde{Z}\in N(k;C)$ with 
$\operatorname{length}Z\leq\operatorname{length}\widetilde{Z}$
such that $\Gamma(\Cal{O}(K_{X}+C))\to 
\Gamma(\Cal{O}(K_{X}+C)\otimes\Cal{O}_{Z})$
is not surjective
but for any subcycle $Z'$ of $Z$ with $\operatorname{length}Z'<
\operatorname{length}Z$,
$\Gamma(\Cal{O}(K_{X}+C))\to 
\Gamma(\Cal{O}(K_{X}+C)\otimes\Cal{O}_{Z'})$
is surjective.
\enddefinition
First we prove the following Theorem.
\proclaim{Theorem 5.2}
Let $X$ be a minimal smooth projective surface
with $\kappa(X)=2$, and let $C$ be an irreducible reduced
curve on $X$ with $C^{2}>0$.
We put $g(C)=q(X)+m$.
We assume that $K_{X}+C$ is not $k$-very ample for some 
integer $k\geq (1/2)(m-1)$.
\flushpar
Assume that 
$$\sharp \bigcup_{\widetilde{Z}\in N(k:C)}
\left(\bigcup_{Z\in S(\widetilde{Z};C)}
\operatorname{Supp}Z\right)=\infty.$$

\comment
Then by Corollary 2.3 in \cite{BeFS}, there exist 
the following {\rm ($\ast$)} and {\rm ($\ast\ast$)};
\roster
\item"($\ast$)"
a 0-dimensional subscheme $Z$ with $\deg Z\leq k+1$ and 
$\operatorname{Supp}(Z)\subset C$;
\item"($\ast\ast$)"
an effective divisor $D_{Z}$ on $X$ such that 
$\operatorname{Supp}(Z)\subset D_{Z}$ and 
$C-2D_{Z}$ is a $\Bbb{Q}$-effective divisor.
\endroster
We also assume $$\sharp \bigcup_{Z\in\Bbb{T}}\operatorname{Supp}(Z)=\infty,$$
where $\Bbb{T}=\{Z\ |\ \text{$Z$ satisfies the above {\rm ($\ast$)}}\}$.
\flushpar
\endcomment

Then $C^{2}\leq 4(k+1)$.
\endproclaim
\demo{Proof}
We remark that $C$ is nef and big.

\comment
First we prove the following Claim.
\proclaim{Claim 5.2}
$K_{X}+C$ is not $m$-spanned at any distinct $m+1$ points on $C$.
\endproclaim
\demo{Proof}
Let $W=\operatorname{Im}(H^{0}(K_{X}+C)\to H^{0}(\omega_{C}))$,
where $\omega_{C}$ is a dualizing sheaf of $C$.
We remark that $\omega_{C}$ is a Cartier divisor.
Let $\alpha$ be the map $H^{0}(K_{X}+C)\to W$.
Then $\dim W=h^{0}(K_{X}+C)-h^{0}(K_{X})=m$
by Riemann-Roch Theorem and Kawamata-Viehweg Vanishing
Theorem.
Let $P_{1},\cdots, P_{m+1}$ be any m+1 distinct points
on $C$.
Let $Z$ be a 0-diminsional subscheme such that
\roster
\item $\Cal{I}_{Z}\Cal{O}_{X,y}=\Cal{O}_{X,y}$ if 
$y\not\in\{ P_{1},\cdots, P_{m+1}\}$;
\item $\Cal{I}_{Z}\Cal{O}_{X,y}=(x_{i},y_{i})$ if 
$y=P_{i}$,
\endroster
where $\Cal{I}_{Z}$ is the ideal sheaf of $Z$ and $(x_{i},y_{i})$
is a local coordinate of $X$ at $P_{i}$.
Let $\beta$ be the restriction map $W\to H^{0}((K_{X}+C)\otimes\Cal{O}_{Z})$.
If $K_{X}+C$ is $m$-spanned at $Z$,
then the restriction 
$\gamma :H^{0}(K_{X}+C)\to H^{0}((K_{X}+C)\otimes\Cal{O}_{Z})$ is surjective.
But we have $\dim W=m$ and $\dim H^{0}((K_{X}+C)\otimes\Cal{O}_{Z})=m+1$.
This is a contradiction since $\gamma=\beta\circ\alpha$.
This completes the proof of Claim 5.2.
\enddemo
We prove Theorem 5.1.
\flushpar
\endcomment

Assume that $C^{2}> 4(k+1)$. 
Then we remark that $C^{2}\geq 2m+3$ by hypothesis.
\flushpar
If $q(X)\leq 2$, then $K_{X}C\geq 0\geq 2q(X)-4$
and so we have $C^{2}\leq 2m+2$ and this is a contradiction.
Hence we have $q(X)\geq 3$.
\flushpar
Then by Corollary 2.3 in \cite{BeFS}, 
for any $Z\in \bigcup\limits_{\widetilde{Z}\in N(k;C)}S(\widetilde{Z};C)$
there exists 
an effective divisor $D_{Z}$ on $X$ such that 
$\operatorname{Supp}(Z)\subset D_{Z}$ and 
$C-2D_{Z}$ is a $\Bbb{Q}$-effective divisor.

\comment
By Corollary 2.3 in \cite{BFS}, there exist 
the following ($\ast$) and ($\ast\ast$);
\roster
\item"($\ast$)"
a 0-dimensional subscheme $Z$ with $\deg Z\leq m+1$ and 
$\operatorname{Supp}(Z)\subset C$;
\item"($\ast\ast$)"
an effective divisor $D_{Z}$ on $X$ such that 
$\operatorname{Supp}(Z)\subset D_{Z}$ and 
$C-2D_{Z}$ is a $\Bbb{Q}$-effective divisor.
\endroster
Then by hypothesis, we have 
$$\sharp \bigcup_{Z\in\Bbb{T}}\operatorname{Supp}(Z)=\infty,$$
where $\Bbb{T}=\{Z\ |\ \text{$Z$ satisfies the above ($\ast$)}\}$.
\flushpar
\endcomment
 
Let
$A=\{ D_{Z}\ |\ \text{$Z\in \bigcup\limits_{\widetilde{Z}\in N(k;C)}
S(\widetilde{Z};C)$ and $D_{Z}$ as above} \}.$
\proclaim{Claim 5.3}
Let $D$ be an effective divisor on $X$ and let $D=\sum\limits_{i}r_{i}C_{i}$
be its prime decomposition.
If there exists an irreducible component $C_{i}$ with $C_{i}^{2}>0$,
and $C-2D$ is $\Bbb{Q}$-effective, then $C^{2}\leq 2m$
if $g(C)=q(X)+m$.
\endproclaim
\demo{Proof}
By Proposition 1.7, we have
$$
\split
K_{X}D&\geq K_{X}C_{i} \\
      &\geq\frac{3}{2}q(X)-3 \\
      &=q(X)+\frac{1}{2}q(X)-3.
\endsplit
$$
Since $q(X)\geq 3$, we obtain that $K_{X}D\geq q(X)-(3/2)$.
Hence $K_{X}D\geq q(X)-1$ because $K_{X}D$ is an integer.
Because $K_{X}$ is nef and $C-2D$ is $\Bbb{Q}$-effective,
we obtain
$$\split
g(C)&=1+\frac{1}{2}(K_{X}+C)C \\
    &\geq 1+\frac{1}{2}(K_{X})(2D)+\frac{1}{2}C^{2} \\
    &=1+K_{X}D+\frac{1}{2}C^{2} \\
    &\geq q(X)+\frac{1}{2}C^{2}.
\endsplit
$$
Therefore $C^{2}\leq 2m$.
This completes the proof of Claim 5.3.
\enddemo
We continue the proof of Theorem 5.2.
\flushpar
By Claim 5.3, any $D_{Z}\in A$ satisfies 
$C_{i}^{2}\leq 0$ for any irreducible component $C_{i}$ of $D_{Z}$.
\flushpar
So $C\not\subset D_{Z}$ for any $D_{Z}\in A$ since $C^{2}>0$.
Hence by hypothesis, we obtain
$$\dim\ \bigcup_{D_{Z}\in A}\left (\bigcup_{C_{Z,i}\in V(D_{Z})}
\operatorname{Supp} C_{Z,i}\right )=2,$$
where $V(D_{Z})=\text{the set of irreducible components of $D_{Z}$}$.
\flushpar
Let
$$ \bigcup_{D_{Z}\in A}V(D_{Z})=B_{1}\cup B_{2},$$
where $B_{1}$ is the set of irreducible curves $C_{1}$ with $C_{1}^{2}<0$
and $B_{2}$ is the set of irreducible curves $C_{2}$ with $C_{2}^{2}=0$.
\flushpar
(1) The case in which $\sharp B_{1}=\infty$.
\flushpar
If $C_{1}\in B_{1}$ with $K_{X}C_{1}\geq q(X)-1$,
then $K_{X}D_{Z}\geq q(X)-1$ and so we have $C^{2}\leq 2m$ by the same
argument as Claim 5.3.
So we have $K_{X}C_{1}\leq q(X)-2$ for any $C_{1}\in B_{1}$.
Then the number of such a curve $C_{1}$ is at most finite by Lemma 1.8.
But this is a contradiction by hypothesis.
\flushpar
(2) The case in which $\sharp B_{2}=\infty$.
\flushpar
If $C_{2}\in B_{2}$ with $K_{X}C_{2}\geq q(X)-1$,
then we have $C^{2}\leq 2m$ by the same
argument as above.
So we have $K_{X}C_{2}\leq q(X)-2$ for any $C_{2}\in B_{2}$.
Then there is a subset $B_{3}\subset B_{2}$ such that 
$\sharp B_{3}=\infty$ and $C_{s}\equiv C_{t}$ for any distinct
$C_{s}, C_{t}\in B_{3}$ by Lemma 1.8.
We take a $C_{k}\in B_{3}$.
Let $\alpha(C_{k})=\dim\operatorname{Ker}(H^{1}(\Cal{O}_{X})\to
H^{1}(\Cal{O}_{C_{k}}))$.
\flushpar
(2-1) The case in which $\alpha(C_{k})\not=0$.
\flushpar
Then by Lemma 1.3 in \cite{Fk4},
there exist an Abelian variety $G$ and
a morphism $f:X\to G$ such that $f(X)$ is 
not a point and $f(C_{k})$ is a point.
Since $C_{k}^{2}=0$, we obtain $f(X)$ is a curve.
By taking Stein factorization, if necessary,
there is a smooth curve $B$, a surjective morphism
$h: X\to B$ with connected fibers,
and a finite morphism $\delta :B\to f(X)$
such that $f=\delta\circ h$.
On the other hand, for any $C_{n}\in B_{3}$ and $C_{n}\not=C_{k}$,
we have $C_{n}C_{k}=C_{k}^{2}=0$.
Hence any element $C_{n}$ of $B_{3}$ is contained in a fiber of $h$
and $C_{n}^{2}=0$.
Therefore for a general fiber $F_{h}$ of $h$, 
we may assume $F_{h}\in B_{3}$.
On the other hand, we have $C-2D_{Z}\leq C-2F_{h}$.
So we obtain that $C-2F_{h}$ is a $\Bbb{Q}$-effective divisor.
\flushpar
Hence we have
$$\split
g(C)&=g(B)+\frac{1}{2}(K_{X/B}+C)C+(CF_{h}-1)(g(B)-1) \\
    &\geq  g(B)+\frac{1}{2}(K_{X/B})(2F_{h})+\frac{1}{2}C^{2} \\
    &=g(B)+2g(F_{h})-2+\frac{1}{2}C^{2} \\
    &=g(B)+g(F_{h})+\frac{1}{2}C^{2}+g(F_{h})-2 \\
    &\geq q(X)+\frac{1}{2}C^{2}
\endsplit
$$
because $K_{X/B}$ is nef, $g(B)\geq 1$ and $g(F_{h})\geq 2$.
\flushpar
Hence $C^{2}\leq 2m$.
But this is a contradiction because we assume that $C^{2}\geq 2m+3$.
\flushpar
(2-2) The case in which $\alpha(C_{k})=0$.
\flushpar
Then $q(X)\leq h^{1}(\Cal{O}_{C_{k}})=g(C_{k})$.
On the other hand, since 
$K_{X}$ is nef, $C_{k}^{2}=0$,
$C-2C_{k}\geq C-2D_{Z}$, and $C-2D_{Z}$
is $\Bbb{Q}$-effective, we obtain
$$
\split
g(C)&=1+\frac{1}{2}(K_{X}+C)C \\
    &\geq 1+\frac{1}{2}(K_{X})(2C_{k})+\frac{1}{2}C^{2} \\
    &=1+K_{X}C_{k}+\frac{1}{2}C^{2} \\
    &=1+2g(C_{k})-2+\frac{1}{2}C^{2} \\
    &\geq 2q(X)-1+\frac{1}{2}C^{2}.
\endsplit
$$
\flushpar
Hence 
$$\split
C^{2}&\leq 2m+2(1-q(X)) \\
     &\leq 2m-4
\endsplit
$$
since $q(X)\geq 3$.
\flushpar
But this is a contradiction by hypothesis.
Therefore $C^{2}\leq 4(k+1)$.
This completes the proof of Theorem 5.2. \qed
\enddemo
\proclaim{Corollary 5.4} Let $X$ be a minimal smooth projective surface with 
$\kappa(X)=2$ and let $C$ be an irreducible reduced curve with $C^{2}>0$.
Then $C^{2}\leq 4m+4$ if $m=g(C)-q(X)$.
\endproclaim
\demo{Proof}
We use Notation 5.1.
By Theorem 5.2, it is sufficient to prove that 
$K_{X}+C$ is not $m$-very ample and 
$$\sharp \bigcup_{\widetilde{Z}\in N(m;C)}\left(\bigcup_{Z\in 
S(\widetilde{Z};C)}\operatorname{Supp}Z\right)=\infty .$$
Let $W=\operatorname{Im}(H^{0}(K_{X}+C)\to H^{0}(\omega_{C}))$,
where $\omega_{C}$ is a dualizing sheaf of $C$.
We remark that $\omega_{C}$ is a Cartier divisor.
Let $\alpha$ be the map $H^{0}(K_{X}+C)\to W$.
Then $\dim W=h^{0}(K_{X}+C)-h^{0}(K_{X})=m$
by Riemann-Roch Theorem and Kawamata-Viehweg Vanishing
Theorem.
Let $P_{1},\cdots, P_{m+1}$ be any m+1 distinct points
on $C\backslash \operatorname{Sing}C$,
where $\operatorname{Sing}C$ denotes the singular locus of $C$.
Let $Z$ be a 0-diminsional subscheme such that
\roster
\item $\Cal{I}_{Z}\Cal{O}_{X,y}=\Cal{O}_{X,y}$ if 
$y\not\in\{ P_{1},\cdots, P_{m+1}\}$;
\item $\Cal{I}_{Z}\Cal{O}_{X,y}=(x_{i},y_{i})$ if 
$y=P_{i}$,
\endroster
where $\Cal{I}_{Z}$ is the ideal sheaf of $Z$ and $(x_{i},y_{i})$
is a local coordinate of $X$ at $P_{i}$
such that $C$ is defined by $(x_{i})$ at $P_{i}$.
Let $\beta$ be the restriction map $W\to H^{0}((K_{X}+C)\otimes\Cal{O}_{Z})$.
If $K_{X}+C$ is $m$-very ample at $Z$,
then the restriction 
$\gamma :H^{0}(K_{X}+C)\to H^{0}((K_{X}+C)\otimes\Cal{O}_{Z})$ is surjective.
But we have $\dim W=m$ and $\dim H^{0}((K_{X}+C)\otimes\Cal{O}_{Z})=m+1$.
This is a contradiction since $\gamma=\beta\circ\alpha$.
Hence $K_{X}+C$ is not $m$-very ample for any 0-dimensional subscheme
with length $m+1$
which consists of distinct $m+1$ points of $C\backslash\operatorname{Sing}(C)$.
This implies 
$$\sharp \bigcup_{\widetilde{Z}\in N(m;C)}\left(\bigcup_{Z\in 
S(\widetilde{Z};C)}\operatorname{Supp}Z\right)=\infty.$$
This completes the proof of Corollary 5.4. \qed
\enddemo
By Corollary 4.27, in order to solve Conjecture 1 (or Conjecture 1$'$),
it is sufficient to consider the case in which $D$ is the following type
($\star$);
\roster
\item"($\star$)"
$D=C_{1}+\sum\limits_{j\geq 2}r_{j}C_{j}$ ;
$C_{1}^{2}>0$ and the intersection matrix 
$\Vert C_{j}, C_{k}\Vert_{j\geq 2, k\geq 2}$
of $\sum\limits_{j\geq 2}r_{j}C_{j}$ is
negative semidefinite.
\endroster
\proclaim{Corollary 5.5}
Let $X$ be a minimal smooth projective surface with
$\kappa(X)=2$ and let $D$ be a nef-big effective divisor on $X$
such that $D$ is the type ($\star$).
Then $D^{2}\leq 4m+4$ if $m=g(D)-q(X)$.
\endproclaim
\demo{Proof}
First we obtain 
$$g(C_{1})=q(X)+m-\frac{1}{2}(K_{X}+D+C_{1})
(\sum\limits_{j\geq 2}r_{j}C_{j}).$$
By Corollary 5.4, we have
$$\split
C_{1}^{2}&\leq 4m+4-2(K_{X}+D+C_{1})
(\sum\limits_{j\geq 2}r_{j}C_{j}) \\
&\leq 4m+4-2(D+C_{1})
(\sum\limits_{j\geq 2}r_{j}C_{j}).
\endsplit
$$
Hence
$$
\split
D^{2}&=C_{1}^{2}+(D+C_{1})(\sum\limits_{j\geq 2}r_{j}C_{j}) \\
&\leq 4m+4-(D+C_{1})(\sum\limits_{j\geq 2}r_{j}C_{j}).
\endsplit
$$
On the other hand
$D+C_{1}$ is nef.
Hence $(D+C_{1})(\sum\limits_{j\geq 2}r_{j}C_{j})\geq 0$
and so we obtain $D^{2}\leq 4m+4$.
This completes the proof of Corollary 5.5. \qed
\enddemo

\subhead{$\S$ 6. Higher dimensional case and Conjecture}
\endsubhead
\par
In this section we consider the case in which $n=\dim X\geq 3$ and
$\kappa(X)\geq 0$.
\proclaim{Theorem 6.1}
Let $(X,L)$ be a quasi-polarized manifold with $\dim X=n\geq 3$
and $\kappa(X)=$0 or 1.
Then $K_{X}L^{n-1}\geq 2(q(X)-n)$.
\endproclaim
\demo{Proof}
\flushpar
(1) The case in which $\kappa(X)=0$.
\flushpar
Then $q(X)\leq n$ by \cite{Ka}.
Hence $K_{X}L^{n-1}\geq 0\geq 2(q(X)-n)$.
\flushpar
(2) The case in which $\kappa(X)=1$.
\flushpar
By Iitaka Theory (\cite{Ii}), 
there exist a smooth projective variety $X_{1}$,
a birational morphism $\mu_{1} :X_{1}\to X$,
a smooth curve $C$,
and a fiber space $f_{1} :X_{1}\to C$ such that $\kappa(F_{1})=0$,
where $F_{1}$ is a general fiber of $f_{1}$.
Let $L_{1}=\mu_{1}^{*}L$.
\flushpar
(2-1) The case in which $g(C)\geq 1$.
\flushpar
By Lemma 1.3.1 and Remark 1.3.2 in \cite{Fk2}
and the semipositivity of $(f_{1})_{*}(K_{X_{1}/C})$ (\cite{Fj1}),
we have $K_{X_{1}/C}L_{1}^{n-1}\geq 0$.
Therefore 
$$
\split
K_{X}L^{n-1}&=K_{X_{1}}L_{1}^{n-1} \\
            &=K_{X_{1}/C}L_{1}^{n-1}+(2g(C)-2)L_{1}^{n-1}F_{1} \\
            &\geq 2g(C)-2.
\endsplit
$$
On the other hand, $q(X)\leq g(C)+(n-1)$ since $q(F_{1})\leq n-1$
by \cite{Ka}.
Hence 
$$\split
K_{X}L^{n-1}&\geq 2(g(C)-1) \\
            &\geq 2(q(X)-n).
\endsplit
$$
\flushpar
(2-2) The case in which $g(C)=0$.
\flushpar
Then $q(X)\leq n-1$ since $q(F_{1})\leq n-1$.
Therefore $K_{X}L^{n-1}\geq 0>2(q(X)-n)$.
\flushpar
This completes the proof of Theorem 6.1. \qed
\enddemo
By considering the above Theorem,
we propose the following Conjecture which is a generalization of 
Conjecture 1$'$.
\proclaim{Conjecture 6.2}
Let $(X,L)$ be a quasi-polarized manifold with $n=\dim X\geq 3$
and $\kappa(X)\geq 0$.
Then $K_{X}L^{n-1}\geq 2(q(X)-n)$.
\endproclaim
By Theorem 6.1, this Conjecture is true if $\kappa(X)=$0 or 1.
We will study Conjecture 6.2 in a future paper.

\enddocument